\documentclass[12pt]{amsart}

\usepackage[german,english]{babel}

\usepackage{dirtytalk}

\usepackage{amsmath}

\usepackage[mathx,matha]{mathabx}

\usepackage{amstext}

\usepackage{amssymb}

\usepackage{array}

\usepackage{amsthm}

\usepackage{dsfont}

\usepackage[latin1]{inputenc}

\usepackage[T1]{fontenc}

\usepackage{mathtools}

\usepackage{appendix}

\usepackage{xcolor}

\usepackage{ulem}

\usepackage[arrow, matrix, curve]{xy}

\usepackage{verbatim}

\usepackage{graphicx}
\graphicspath{{images/}}

\usepackage{subfigure}

\usepackage{hyperref}

\usepackage{cases}

\usepackage[update,prepend]{epstopdf}

\usepackage{bbold}

\usepackage[a]{esvect}

\usepackage{tikz}

\newtheorem{theorem}{Theorem}[section]

\newtheorem*{theorem*}{Theorem}

\theoremstyle{plain}

\newtheorem{conjecture}[theorem]{Conjecture}

\newtheorem*{conjecture*}{Conjecture}

\newtheorem{assumption}[theorem]{Assumption}

\newtheorem{corollary}[theorem]{Corollary}

\newtheorem{proposition}[theorem]{Proposition}

\newtheorem{lemma}[theorem]{Lemma}

\theoremstyle{remark}

\newtheorem{condition}{Condition}

\theoremstyle{definition}

\newtheorem{definition}[theorem]{Definition}

\newtheorem{remark}[theorem]{Remark}

\newcommand{\orb}{{\rm orb}}

\def\eps{\varepsilon}

\def\bi{\begin{itemize}}

\def\ei{\end{itemize}}

\newcommand{\R}{\mathbb{R}}

\newcommand{\N}{\mathbb{N}}

\newcommand{\Z}{\mathbb{Z}}

\newcommand{\T}{\mathbb{T}}

\newcommand{\A}{\mathbb{A}}

\renewcommand{\phi}{\varphi}

\DeclareMathOperator{\vol}{vol}

\DeclareMathOperator{\conv}{conv}

\DeclareOldFontCommand{\it}{\normalfont\itshape}{\mathit}

\newcommand{\bspm}{\left(\begin{smallmatrix}}\newcommand{\espm}{\end{smallmatrix}\right)}

\newcommand{\bpm}{\begin{pmatrix}}\newcommand{\epm}{\end{pmatrix}}

\def\bs{\begin{satz}}\def\es{\end{satz}}

\def\blem{\begin{lemma}}\def\elem{\end{lemma}}

\def\bthm{\begin{theorem}}\def\ethm{\end{theorem}}

\def\bcor{\begin{corollary}}\def\ecor{\end{corollary}}

\def\beq{\begin{equation}}\def\eeq{\end{equation}}

\def\beqq{\begin{equation*}}\def\eeqq{\end{equation*}}

\def\bal{\begin{align}}\def\eal{\end{align}}

\def\bpf{\begin{proof}}\def\epf{\end{proof}}

\def\bex{\begin{example}}\def\eex{\end{example}}

\def\brem{\begin{remark}}\def\erem{\end{remark}}

\def\bass{\begin{assumption}}\def\eass{\end{assumption}}

\def\bprop{\begin{proposition}}\def\eprop{\end{proposition}}

\def\bdefi{\begin{definition}}\def\edefi{\end{definition}}

\def\bcond{\begin{condition}}\def\econd{\end{condition}}

\def\bconj{\begin{conjecture}}\def\econj{\end{conjecture}}

\DeclareFontEncoding{FMS}{}{}

\DeclareFontSubstitution{FMS}{futm}{m}{n}

\DeclareFontEncoding{FMX}{}{}

\DeclareFontSubstitution{FMX}{futm}{m}{n}

\DeclareSymbolFont{fouriersymbols}{FMS}{futm}{m}{n}

\DeclareSymbolFont{fourierlargesymbols}{FMX}{futm}{m}{n}

\DeclareMathDelimiter{\VERT}{\mathord}{fouriersymbols}{152}{fourierlargesymbols}{147}

\def\bi{\begin{itemize}}

\def\ei{\end{itemize}}

\def\ben{\begin{enumerate}}

\def\een{\end{enumerate}}

\usepackage{enumitem}

\usepackage{tcolorbox}

\newtcolorbox{implementation}[2][]{colframe=blue!75!black,colbacktitle=green!10!white,colback=green!10!white,coltitle=green!75!black,title={#2},fonttitle=\bfseries,#1}

\begin{document}

\title[Viterbo's conjecture as a worm problem]{Viterbo's conjecture as a worm problem}

\author{Daniel Rudolf}


\date{\today}

\maketitle

\begin{abstract}
In this paper, we relate Viterbo's conjecture from symplectic geometry to Minkowski versions of worm problems which are inspired by the well-known Moser worm problem from geometry. For the special case of Lagrangian products this relation provides a connection to systolic Minkowski billiard inequalities and Mahler's conjecture from convex geometry. Moreover, we use the above relation in order to transfer Viterbo's conjecture to a conjecture for the longstanding open Wetzel problem which also can be expressed as a systolic Euclidean billiard inequality and for which we discuss an algorithmic approach in order to find a new lower bound. Finally, we point out that the above mentioned relation between Viterbo's conjecture and Minkowski worm problems has a structural similarity to the known relationship between Bellmann's lost-in-a-forest problem and the original Moser worm problem.
\end{abstract}

\section{Introduction and main results}\label{Sec:Introworm}

Worm problems have a long history. The earliest known problem of this type was posed by L. Moser in \cite{Moser1966} (cf.\;also \cite{Moser1991}) more than 50 years ago:\\

\par
\begingroup
\leftskip=1cm
\noindent \textit{Moser's worm problem: Find a/the (convex) set of least area that contains a congruent copy of each arc in the plane of lenth one.}\\
\par
\endgroup

\noindent Here, the unit arcs are sometimes called worms, while the problem has been phrased in many different ways in the literature: the architect's version (find the smallest comfortable living quarters for a unit worm), the humanitarian version (find the shape of the most efficient worm blanket), the sadistic version (find the shape of the best mallet head), and so on (cf.\;\cite{Wetzel2003}). So far, despite a lot of research, only partial results are known, including the existence of such a minimum cover in the convex case (probably the first time proven in \cite{LaiPoo1986}), but its shape and area remain unknown. The best bounds presently known for its area $\mu$ are:
\beqq 0.23224\leq \mu \leq 0.27091\footnote{We round all decimal numbers up to the fifth decimal place.}\eeqq
(cf.\;\cite{KhPagSri2013} for the lower and \cite{Wang2006} for the upper bound).

Worm problems can be formulated in considerable generality (cf.\;\cite{Wetzel2003}):\\

\par
\begingroup
\leftskip=1cm 
\noindent\textit{Given a collection $\mathcal{F}$ of $n$-dimensional figures $F$ and a transitive group $\mathcal{M}$ of motions $m$ on $\R^n$, find minimal convex target sets $K\subset \R^n$--minimal in the sense of having least volume, surface volume, or whatever--so that for each $F\in\mathcal{F}$ there is a motion $m \in \mathcal{M}$ with
\beqq m(F)\subseteq K.\eeqq}
\par
\endgroup

The existence of solutions to this problem can be guaranteed under certain natural hypotheses by fundamental compactness results like the Blaschke selection theorem (cf.\;\cite[\S 18]{Blaschke1916} for Blaschke's selection theorem and \cite{KellyWeiss1979} or \cite{LaiPoo1986} for its application; cf.\;also Theorem \ref{Thm:Blaschkeselection} and its application in Propositions \ref{Prop:minimumexists}, \ref{Prop:minimaxproblem1}, \ref{Prop:minimaxproblem2}, and \ref{Prop:minimaxproblem3}).

When the problem does not permit an arc to be replaced by its mirror image, then it is appropriate to consider the subgroup of orientation preserving motions. For other problems, e.g., Moser's original worm problem, orientation reversing motions are permitted. Many problems whose motion group is the group of translations have been studied in the literature (cf.\;\cite{BezCon1989}, \cite{CroftFalconerGuy1991}, \cite{Wetzel1973}).

In order to formulate the specific worm problem which is of main interest for our study, we introduce the following definition: Let $T\subset\R^n$ be a convex body, i.e., a compact convex set in $\R^n$ with nonempty interior, and $T^\circ$ its polar. Using the Minkowski functional $\mu_{T^\circ}$ with respect to $T$'s polar, we define the \textit{$\ell_T$-length} of a closed $H^1([0,\widetilde{T}],\R^n)$-curve\footnote{This implies that $q$ is differentiable almost everywhere with $\dot{p}\in L^2([0,\widetilde{T}],\R^n)$.} $q$ (from now on, for the sake of simplicity, every closed curve is assumed to fulfill this Sobolev property), $\widetilde{T}\geq 0$, by
\beqq \ell_T(q):=\int_0^{\widetilde{T}} \mu_{T^\circ}(\dot{p}(t))\, \mathrm{d}t.\eeqq

The worm problem which is of main interest for our study we call the \textit{Minkowski worm problem}. Referring to the above general worm problem formulation, for this for convex body $T\subset\R^n$, we consider $\mathcal{F}=\mathcal{F}(T,\alpha)$ as the set of closed curves of $\ell_T$-length $\alpha > 0$, $\mathcal{M}$ as the group of translations and the minimization in the sense of having minimal volume:\\

\par
\begingroup
\leftskip=1cm 
\noindent \textit{Minkowski worm problem: Let $T\subset\R^n$ be a convex body. Find the volume-minimizing convex bodies $K\subset\R^n$ that contain a translate of every closed curve of $\ell_T$-length $\alpha$.}\\
\par
\endgroup

\noindent So, in contrast to Moser's worm problem, we consider general dimension (instead of just dimension two), length-measuring with Minkowski functionals with respect to arbitrary convex bodies (instead of with respect to the Euclidean unit ball), closed curves (instead of not necessarily closed arcs), and translations (instead of congruence transformations). In other words and introducing a notation which will be useful throughout this paper: Let $cc(\R^n)$ be the set of closed curves in $\R^n$. Find the minimizers\footnote{In Proposition \ref{Prop:minimumexists}, we will prove that in fact there exists at least one minimizer.} of
\beqq \min_{K\in A(T,\alpha)} \vol(K),\quad \eeqq
where for convex body $T\subset\R^n$ and $\alpha > 0$, we define
\beqq A(T,\alpha):=\left\{K\subset\R^n \text{ convex body}: L_T(\alpha)\subseteq C(K)\right\}\eeqq
with
\beqq L_T(\alpha):=\left\{q \in cc(\R^n) : \ell_T(q)=\alpha\right\}\eeqq
and
\beqq C(K):=\left\{q \in cc(\R^n) : \exists k\in\R^n \text{ s.t. } q \subseteq k + K\right\},\eeqq
where, for the sake of simplicity, we, in general, identify $q$ with its image.

\begin{figure}[h!]
\centering
\def\svgwidth{250pt}
\includegraphics[scale=.9]{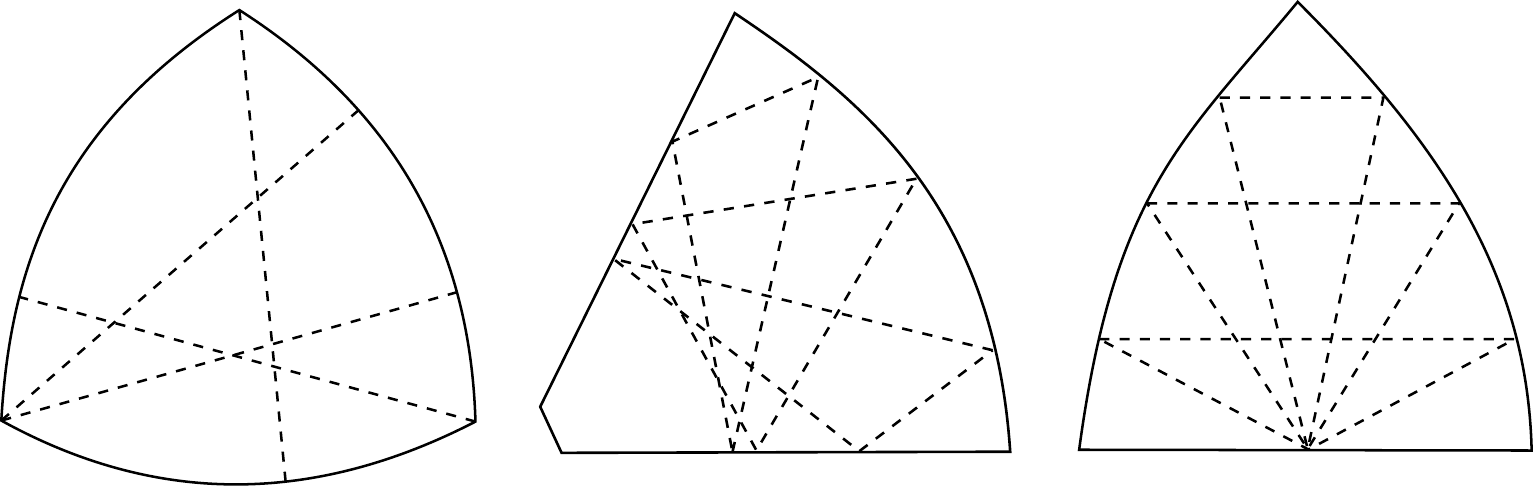}
\caption[Three examples of convex bodies approaching the minimizer in Wetzel's problem]{On the left side is the Reuleaux triangle with width $\frac{1}{2}$ and volume $0.17619$, in the middle is a convex body with volume $0.17141$ which was found by Wetzel in \cite{Wetzel1973}, and on the right is a convex body, looking a bit like a church window, with base length and height equal to $\frac{1}{2}$ and volume $\frac{1}{6}\approx 0.16667$ (for both the middle and right convex body we refer to \cite{BezCon1989}). Some worms are drawn in in each case.}
\label{img:upperbounds}
\end{figure}

The case when the dimension is $2$, $T$ is the Euclidean unit ball in $\R^2$, and, without loss of generality, $\alpha =1$ is the only case (one could say: the two-dimensional Euclidean worm problem) which has been investigated so far. It is known as \textit{Wetzel's problem}. So far, the minimal volume (area) for this problem is not known, but the best bounds presently known for the minimum are $0.15544$ as lower (cf.\;\cite{Wetzel1973}, where an argument from \cite{Pal1921} is used) and $0.16526$ as upper bound (cf.\;\cite{BezCon1989}; note that in \cite{Wetzel1973} it was claimed incorrectly an upper bound of $0.159$). In comparison to that: The volumes of the obvious covers of constant width, the ball of radius $1/4$ and the Reuleaux triangle of width $1/2$, are $0.19635$ and $0.17619$, respectively. Since, by the Blaschke-Lebesgue theorem, the Reuleaux triangle is the volume-minimizing set of constant width (cf.\;\cite{Blaschke1915} and \cite{Lebesgue1914}; cf.\;\cite{Harrell2002} for a direct proof by analyzing the underlying variational problem), we can conclude that a minimizer for Wetzel's problem is not of constant width. We refer to Figure \ref{img:upperbounds} for three examples whose volumes are approaching (not achieving) the minimum (clearly, the middle and right convex bodies are not of constant width).

Although we derive some results, the primary goal of our study will not be to solve these Minkowski worm problems, rather to relate them to Viterbo's conjecture from symplectic geometry (cf.\;\cite{Viterbo2000}) which for convex bodies $C\subset\R^{2n}$ reads
\beqq \vol(C)\geq \frac{c_{EHZ}(C)^n}{n!}.\eeqq
For that, we recall that the EHZ-capacity of a convex body $C\subset\R^{2n}$ can be defined\footnote{This definition is in fact the outcome of a historically grown study of symplectic capacities. Traced back--recalling that $c_{EHZ}$ in its present form is the generalization of a symplectic capacity by K\"{u}nzle in \cite{Kuenzle1996} after applying the dual action functional introduced by Clarke in \cite{Clarke1979}--, the EHZ-capacity denotes the coincidence of the Ekeland-Hofer- and Hofer-Zehnder-capacities, originally constructed in \cite{EkHo1989} and \cite{HoferZehnder1990}, respectively.} by
\beqq c_{EHZ}(C)=\min\{\A(x):x \text{ closed characteristic on }\partial C\},\eeqq
where a closed characteristic on $\partial C$ is an absolutely continuous loop in $\R^{2n}$ satisfying
\beqq \begin{cases}\dot{x}(t)\in J\partial H_{C}(x(t))\quad \text{a.e}. \\ H_{C}(x(t))=\frac{1}{2}\;\forall t\in \T\end{cases}\eeqq
where
\beqq H_{C}(x)=\frac{1}{2}\mu_{C}(x)^2,\quad J=\begin{pmatrix} 0 & \mathbb{1} \\ -\mathbb{1} & 0 \end{pmatrix},\quad \T=\R/\widetilde{T}\Z,\; \widetilde{T}>0.\eeqq
$\widetilde{T}$ is the period of the loop and by $\A$ we denote its action defined by
\beqq \A(x)= -\frac{1}{2}\int_0^{\widetilde{T}} \langle J\dot{x}(t),x(t)\rangle \, \mathrm{d}t.\eeqq

The first main result of this paper addresses the special case of Lagrangian products
\beqq C=K\times T\subset \R^n_q \times \R^n_p \cong \R^{2n},\eeqq
where $K$ and $T$ are convex bodies in $\R^n$. We denote by $\mathcal{C}(\R^n)$ the set of convex bodies in $\R^n$.

\bthm\label{Thm:main1}
Viterbo's conjecture for convex Lagrangian products $K\times T\subset\R^n \times \R^n$
\beqq \vol(K\times T) \geq \frac{c_{EHZ}(K\times T)^n}{n!},\quad K,T \in \mathcal{C}(\R^n),\eeqq
is equivalent to
\beq \min_{\vol(T)=1}\;\min_{K\in A(T,1)}\; \vol(K) \geq \frac{1}{n!},\quad K,T \in \mathcal{C}(\R^n).\footnote{In order not to let the entries under the maxima and minima become too long, we decided in general to write extremization problems of the form
\beqq \min_{\vol(T)=1,\; T\in \mathcal{C}(\R^n)}\;\min_{K\in A(T,1),\, K\in \mathcal{C}(\R^n)}\; \vol(K)\eeqq
(here, $K\in\mathcal{C}(\R^n)$ is already included within $K\in A(T,1)$; however, it happens sometimes that only convex and centrally symmetric convex bodies are considered in which case the formulations of the extremization problems become as long as indicated) as
\beqq \min_{\vol(T)=1}\;\min_{K\in A(T,1)}\; \vol(K),\quad K,T\in \mathcal{C}(\R^n).\eeqq}\label{eq:main11}\eeq
Additionally, equality cases $K^*\times T^*$ of Viterbo's conjecture satisfying
\beqq \vol(K^*)=\vol(T^*)=1\eeqq
are composed of equality cases $(K^*,T^*)$ of \eqref{eq:main11}. Conversely, equality cases $(K^*,T^*)$ of \eqref{eq:main11} form equality cases $K^*\times T^*$ of Viterbo's conjecture.
\ethm

This yields the following corollary, which seems to be more suitable in order to approach Viterbo's conjecture as an optimization problem (cf.\;Section \ref{Sec:OptimizationProblem}).

\bcor\label{Cor:optimizationmain1}
Viterbo's conjecture for convex Lagrangian products $K\times T\subset\R^n \times \R^n$
\beqq \vol(K\times T) \geq \frac{c_{EHZ}(K\times T)^n}{n!},\quad K,T \in \mathcal{C}(\R^n),\eeqq
is equivalent to
\beq \min_{\vol(T)=1}\;\min_{a_q\in\R^n}\vol\bigg(\conv\bigg\{ \bigcup_{q\in L_T(1)} (q+a_q)\bigg\}\bigg)\geq \frac{1}{n!},\quad T \in \mathcal{C}(\R^n),\footnote{Here, we note that $K$ has been dissolved by replacing it by an expression that extremizes over all possible $K$s. The extremizing $K$ is of the form \eqref{eq:Kstern}.}\label{eq:optimizationmain11}\eeq
where the second minimum on the left runs for every $q\in L_T(1)$ over all possible translations in $\R^n$. Additionally, equality cases $K^*\times T^*$ of Viterbo's conjecture satisfying
\beqq \vol(K^*)=\vol(T^*)=1\eeqq
are composed of equality cases $T^*$ of \eqref{eq:optimizationmain11} with
\beq K^*=\conv\bigg\{ \bigcup_{q\in L_{T^*}(1)} (q+a_q^*) \bigg\},\label{eq:Kstern}\eeq
where $a_q^*$ are the minimizers in \eqref{eq:optimizationmain11}. Conversely, equality cases $T^*$ of \eqref{eq:optimizationmain11} with $K^*$ as in \eqref{eq:Kstern} form equality cases $K^*\times T^*$ of Viterbo's conjecture.
\ecor

In analogy to Theorem \ref{Thm:main1}, also Mahler's conjecture from convex geometry (cf.\;\cite{Mahler1939}), i.e.,
\beq \vol(T)\vol(T^\circ)\geq \frac{4^n}{n!},\quad T\in \mathcal{C}^{cs}(\R^n),\label{eq:mahlerdefinition}\eeq
where by $\mathcal{C}^{cs}(\R^n)$ we denote the set of all centrally symmetric convex bodies in $\R^n$, can be expressed as a worm problem. As shown in \cite{ArtKarOst2013}, this is due to the fact that Mahler's conjecture is a special case of Viterbo's conjecture.

\bthm\label{Thm:Mahler}
Mahler's conjecture for centrally symmetric convex bodies
\beq \vol(T)\vol(T^\circ)\geq \frac{4^n}{n!},\quad T\in \mathcal{C}^{cs}(\R^n),\label{eq:Mahler0}\eeq
is equivalent to
\beq \min_{T\in A\left(T^\circ,\sqrt[n]{\vol(T^\circ)}\right)} \vol(T) \geq \frac{1}{n!},\quad T\in \mathcal{C}^{cs}(\R^n).\label{eq:Mahler1}\eeq
Additionally, equality cases $T^*$ of Mahler's conjecture \eqref{eq:Mahler0} satisfying
\beqq \vol(T^*)=1\eeqq
are equality cases of \eqref{eq:Mahler1}. And conversely, equality cases $T^*$ in \eqref{eq:Mahler1} are equality cases of Mahler's conjecture \eqref{eq:Mahler0}.
\ethm

Furthermore, also \textit{systolic Minkowski billiard inequalities} within the field of billiard dynamics can be related to worm problems. 

In order to state this, let us recall some relevant notions from the theory of Minkowski billiards (cf.\;\cite{KruppRudolf2022}): For convex bodies $K,T\subset\R^n$, we say that a closed polygonal curve\footnote{For the sake of simplicity, whenever we talk of the vertices $q_1,...,q_m$ of a closed polygonal curve, we assume that they satisfy $q_j\neq q_{j+1}$ and $q_j$ is not contained in the line segment connecting $q_{j-1}$ and $q_{j+1}$ for all $j\in\{1,...,m\}$. Furthermore, whenever we settle indices $1,...,m$, then the indices in $\Z$ will be considered as indices modulo $m$.\label{foot:polygonalline}} with vertices $q_1,...,q_m$, $m\geq 2$, on the boundary of $K$ is a \textit{closed weak $(K,T)$-Minkowski billiard trajectory} if for every $j\in \{1,...,m\}$, there is a $K$-supporting hyperplane $H_j$ through $q_j$ such that $q_j$ minimizes
\beqq \mu_{T^\circ}(\widebar{q}_j-q_{j-1})+\mu_{T^\circ}(q_{j+1}-\widebar{q}_j),\eeqq
over all $\widebar{q}_j\in H_j$. We encode this closed $(K,T)$-Minkowski billiard trajectory by $(q_1,...,q_m)$. Furthermore, we say that a closed polygonal curve with vertices $q_1,...,q_m$, $m\geq 2$, on the boundary of $K$ is a \textit{closed (strong) $(K,T)$-Minkowski billiard trajectory} if there are points $p_1,...,p_m$ on $\partial T$ such that
\beqq \begin{cases}q_{j+1}-q_j\in N_T(p_j),\\ p_{j+1}-p_j=-N_K(q_{j+1})\end{cases}\eeqq
is fulfilled for all $j\in\{1,...,m\}$. We denote by $M_{n+1}(K,T)$ the set of closed $(K,T)$-Minkowski billiard trajectories with at most $n+1$ bouncing points.

Then, for convex body $K\subset\R^n$, introducing $F^{cp}(K)$ as the set of all closed polygonal curves in $\R^n$ that cannot be translated into $\mathring{K}$, we have the following relations:

\bthm\label{Thm:relations}
Let $T\subset\R^n$ be a convex body and $\alpha,c >0$. Then, the following statements are equivalent:
\begin{itemize}
\item[(i)]
\beqq \max_{\vol(K)=c}\; \min_{q\in F^{cp}(K)}\ell_T(q) \leq \alpha, \quad K\in\mathcal{C}(\R^n),\eeqq
\item[(ii)]
\beqq \max_{\vol(K)=c}\; c_{EHZ}(K\times T) \leq \alpha, \quad K\in\mathcal{C}(\R^n),\eeqq
\item[(iii)]
\beqq \max_{\vol(K)=c}\;\; \min_{q \in M_{n+1}(K,T)} \ell_{T}(q) \leq \alpha, \quad K\in\mathcal{C}(\R^n),\eeqq
\item[(iv)]
\beqq \min_{K\in A(T,\alpha)} \vol(K)\geq c, \quad K\in\mathcal{C}(\R^n),\eeqq
\item[(v)] 
\beqq \min_{a_q\in\R^n} \vol\bigg(\conv \bigg\{ \bigcup_{q\in L_T(1)}(q+a_q) \bigg\}\bigg) \geq c, \quad K\in\mathcal{C}(\R^n).\eeqq
\end{itemize}
If $T$ is additionally assumed to be strictly convex, then the following systolic weak Minkowski billiard inequality can be added to the above list of equivalent expressions:
\begin{itemize}
\item[(vi)] 
\beqq \max_{\vol(K)=c}\;\; \min_{q \text{ cl.\,weak }(K,T)\text{-Mink.\,bill.\,traj.}} \ell_{T}(q) \leq \alpha, \quad K\in\mathcal{C}(\R^n).\eeqq
\end{itemize}

Moreover, every equality case $(K^*,T^*)$ of any of the above inequalities is also an equality case of all the others.
\ethm

Now, we turn our attention to the general Viterbo conjecture for convex bodies in $\R^{2n}$. For that, we first introduce the following definitions: We denote by $\mathcal{C}^p\left(\R^{2n}\right)$ the set of convex polytopes in $\R^{2n}$. For $P\in\mathcal{C}^p\left(\R^{2n}\right)$, we denote by
\beqq F^{cp}_*(P)\subset F^{cp}(P)\eeqq
the set of all closed polygonal curves $q=(q_1,...,q_m)$ in $F^{cp}(P)$ for which $q_j$ and $q_{j+1}$ are on neighbouring facets $F_j$ and $F_{j+1}$ of $P$ such that there are $\lambda_j,\mu_{j+1}\geq 0$ with
\beqq q_{j+1}=q_j + \lambda_j J\nabla H_P(x_j)+\mu_{j+1}J\nabla H_P(x_{j+1}),\eeqq
where $x_j$ and $x_{j+1}$ are arbitrarily chosen interior points of $F_j$ and $F_{j+1}$, respectively. Later, we will see that the existence of such closed polygonal curves is guaranteed.

\bthm\label{Thm:main2}
Viterbo's conjecture for convex polytopes in $\R^{2n}$
\beq \vol(P)\geq \frac{c_{EHZ}(P)^n}{n!},\quad P \in \mathcal{C}^p\left(\R^{2n}\right),\label{eq:main20}\eeq
is equivalent to
\beq \min_{P\in A\left(JP,\frac{1}{R_{P}}\right)} \vol(P) \geq \frac{1}{2^n n!},\quad P \in \mathcal{C}^p\left(\R^{2n}\right),\label{eq:main21}\eeq
where we define
\beqq R_{P}:=\frac{\min_{q\in F^{cp}_*(P)}\ell_{\frac{JP}{2}}(q)}{\min_{q\in F^{cp}(P)}\ell_{\frac{JP}{2}}(q)} \geq 1.\eeqq
Additionally, $P^*$ is an equality case of Viterbo's conjecture for convex polytopes \eqref{eq:main20} satisfying
\beqq \vol(P^*)=1\eeqq
if and only if $P^*$ is an equality case of \eqref{eq:main21}.
\ethm

When we look at the operator norm of the complex structure/symplectic matrix $J$ with respect to a convex body $C\subset\R^{2n}$ as map from
\beqq \left(\R^{2n},||\cdot||_{C^\circ}\right) \; \text{ to } \; \left(\R^{2n},||\cdot||_{C}\right)\eeqq
as it has been done in \cite{AkopKar2017} and \cite{GlusOst2015}, i.e.,
\beqq ||J||_{C^\circ \rightarrow C}=\sup_{||v||_{C^\circ}\leq 1}||Jv||_C,\eeqq
then we derive the following theorem:

\bthm\label{Thm:main3}
Viterbo's conjecture for convex bodies in $\R^{2n}$
\beq \vol(C) \geq \frac{c_{EHZ}(C)^n}{n!},\quad C\in \mathcal{C}\left(\R^{2n}\right),\label{eq:main300}\eeq
is equivalent to
\beq \min_{C\in A\left(C^\circ,\frac{1}{\widetilde{R}_C}\right)} \vol(C) \geq \frac{1}{n!},\quad C\in \mathcal{C}\left(\R^{2n}\right),\label{eq:main3}\eeq
where
\beqq \widetilde{R}_C:=\frac{c_{EHZ}(C)}{c_{EHZ}(C\times C^\circ)}\geq \frac{1}{2||J||_{C^\circ \rightarrow C}}.\eeqq
Additionally, $C^*$  is an equality case of Viterbo's conjecture for convex bodies in $\R^{2n}$ \eqref{eq:main300} satisfying
\beqq \vol(C^*)=1\eeqq
if and only if $C^*$ is an equality case of \eqref{eq:main3}.
\ethm

Finally, we turn to Wetzel's problem. For that, we keep the current state of things in mind:

\bthm[Wetzel in \cite{Wetzel1973}, '73; Bezdek \& Connelly in \cite{BezCon1989}, '89]\label{Thm:WetzelBezdekConnelly}
In dimension $n=2$, we have
\beqq \min_{K\in A\left(B_1^2,1\right)}\vol(K)\in (0.15544,0.16526),\quad K\in\mathcal{C}(\R^2),\eeqq
where we denote by $B_1^2$ the Euclidean unit ball in $\R^2$.
\ethm

Then, as application of Theorem \ref{Thm:main1}, we transfer Viterbo's conjecture onto Wetzel's problem. This results in the following conjecture:

\bconj\label{Conj:Wetzelsproblem}
We have
\beqq \min_{K\in A\left(B_1^2,1\right)}\vol(K)\geq \frac{1}{2\pi} \approx 0.15915, \quad K\in\mathcal{C}(\R^2).\eeqq
\econj

Applying \cite[Theorem 3.12]{KruppRudolf2022} and Theorem \ref{Thm:relations}, we note that this conjecture can be equivalently expressed as \textit{systolic Euclidean billiard inequality}:

\bconj\label{Conj:WetzelsproblemBilliard}
We have
\beqq \max_{\vol(K)=\frac{1}{2\pi}}\;\; \min_{q \text{ cl.\,}(K,B_1^2)\text{-Mink.\,bill.\,traj.}} \ell_{B_1^2}(q)\leq 1\eeqq
for $K\in \mathcal{C}(\R^2)$.
\econj

We remark that, for the configuration $(K,B_1^2)$, due to the strict convexity of $B_1^2$, the notions of weak and strong $(K,B_1^2)$-Minkowski billiards coincide and are equal to the one of billiards in the Euclidean sense.

Although much work has been done around Wetzel's problem and the systolic Euclidean billiard inequality, this shows that Viterbo's conjecture is even unsolved for the \say{trivial} configuration
\beqq K \times B_1^2\subset\R^2\times\R^2.\eeqq
On the other hand, looking at these two problems from the symplectic point of view, can help us to conceptualize them from a very different point of view.

The worm problems in Theorems \ref{Thm:Mahler}, \ref{Thm:main2} and \ref{Thm:main3} seem to be very hard to solve (as it is expected from the perspective of Mahler's/Viterbo's conjecture). On the one hand, this is a consequence of the inner dependencies within
\beqq T\in A\left(T^\circ,\sqrt[n]{\vol(T^\circ)}\right) ,\; P\in A\left(JP,\frac{1}{R_{P}}\right)  \; \text{ and }\; C\in A\left(C^\circ,\frac{1}{\widetilde{R}_C}\right),\eeqq
on the other hand, the estimates, hidden in $R_P$ and $\widetilde{R}_C$, beyond specific configurations, do not seem to be so accessible. Nevertheless, beyond the inner dependencies of the second arguments of $A(\cdot,\cdot)$, perhaps it turns out to be fruitful to investigate worm problems of the following structure a little bit more in detail: For $\alpha > 0$, find
\beqq \min_{C\in A(C^\circ,\alpha)} \vol(C) \; \text{ and } \; \min_{C\in A(J^{-1}C,\alpha)}\vol(C).\eeqq
Interestingly enough, from this perspective, Viterbo's and Mahler's conjecture are very similar in structure.

Motivated by a relationship between Moser's worm problem and a version of \textit{Bellman's lost-in-a-forest problem} shown by Finch and Wetzel in \cite{FinchWetzel2004}, we further investigate whether it is possible also to relate Minkowski worm problems to versions of Bellman's lost-in-a-forest problem. And indeed, it will turn out that the relationship established in \cite{FinchWetzel2004} is somewhat similar to the relationship between Minkowski worm problems and Viterbo's conjecture for convex Lagrangian products. However, before we will elaborate on this, we will give a short introduction to Bellman's lost-in-a-forest problem and general escape problems of this type.

In 1955, Bellman stated in \cite{Bellman1956} the following research problem (cf.\;also \cite{Bellman1957} and \cite{Bellman1963}):\\

\par
\begingroup
\leftskip=1cm
\noindent \textit{We are given a region $R$ and a random point $P$ within the region. Determine the paths which (a) minimize the expected time to reach the boundary, or (b) minimize the maximum time required to reach the boundary.}\\
\par
\endgroup

\noindent This problem can be phrased as:\\

\par
\begingroup
\leftskip=1cm
\noindent \textit{A hiker is lost in a forest whose shape and dimensions are precisely known to him. What is the best path for him to follow to escape from the forest?}\\
\par
\endgroup

\noindent In other words: To solve the lost-in-a-forest problem one has to find the \textit{best} escape path--the best in terms of minimizing the maximum or expected time required to escape the forest. A third interpretation of \textit{best} has been given in \cite{CroftFalconerGuy1991}: Find the best escape path in terms of maximizing the probability of escape within a specified time period.

Bellman asked about two configurations in particular: on the one hand, the configuration in which the region is the infinite strip between two parallel lines a known distance apart, on the other hand, the configuration in which the region is a half-plane and the hiker's distance from the boundary is known. For the case when \textit{best} is understood in terms of the maximum time to escape, both of these two configurations have been studied: for the first configuration, the best path was found in \cite{Zalgaller1961} ('61), for the second, in \cite{Isbell1957} ('57) (where a complete and detailed proof was not published until it was done in \cite{Joris1980} ('80); cf.\;\cite{Finch2019} for an english translation). In each of these two cases, the shortest escape path is unique up to congruence. Apart from that, not much is known for other interpretations of \textit{best}. We refer to \cite{Ward2008} for a detailed survey on the different types, results, and some related material.

Finch and Wetzel studied in \cite{FinchWetzel2004} the case in which the \textit{best} escape path is the shortest. As already mentioned above, in this case, they could show a fundamental relation to Moser's worm problem.

Before we further elaborate on this, it is worth mentioning to note that Williams in \cite{Williams2002} has included lost-in-a-forest problems in his recent list \say{Million Buck Problems} of unsolved problems of high potential impact on mathematics. He justified the selection of these problems by mentioning that the techniques involved in their resolution will be worth at least one million dollars to mathematics.

Now, let's consider the case studied by Finch and Wetzel and take it a little more rigorously. For that, let $\gamma$ be a \textit{path} in $\R^2$, i.e., a continuous and rectifiable mapping of $[0,1]$ into $\R^2$. Let $\ell_{B_1^2}(\gamma)$ be its Euclidean length and $\{\gamma\}$ its trace $\gamma([0,1])$. We call a \textit{forest} a closed, convex region in the plane with nonempty interior. A path $\gamma$ is an \textit{escape path} for a forest $K$ if a congruent copy of it meets the boundary $\partial K$ no matter how it is placed with its initial point in $K$, i.e., for each point $P\in K$ and each Euclidean motion (translation, rotation, reflection and combinations of them) $\mu$ for which $P=\mu(\gamma(0))$ the intersection $\mu\left(\{\gamma\}\right)\cap \partial K$ is nonempty. Then, among all the escape paths for a forest $K$, there is at least one whose length is the shortest. The \textit{escape length} $\alpha$ of a forest $K$ is the length of one of these shortest escape paths for $K$. Based on these notions, Finch and Wetzel proved the following:

\bthm[Theorem 3 in \cite{FinchWetzel2004}]\label{Thm:FinchWetzel}
Let $K\subset\R^2$ be a convex body. The escape length $\alpha^*$ of $K$ is the largest $\alpha$ for which for every path $\gamma$ with length $\leq \alpha$, there is a Euclidean motion $\mu$ such that $K$ covers $\mu\left(\{\gamma\}\right)$.
\ethm

For Finch and Wetzel, this theorem established the connection to Moser's worm problem. For that, we recall that in Moser's worm problem one tries to find a/the convex set of least area that contains a congruent copy of each arc in the plane of a certain length. Clearly, the condition of having a certain length can be replaced by the condition of having a length which is bounded from above by that certain length.

Now, translated into our setting, we can derive a similar result. For that, we first have to define a version of a lost-in-a-forest problem which is compatible with the Minkowski worm problems discussed in the previous sections.

\begin{figure}[h!]
\centering
\def\svgwidth{375pt}
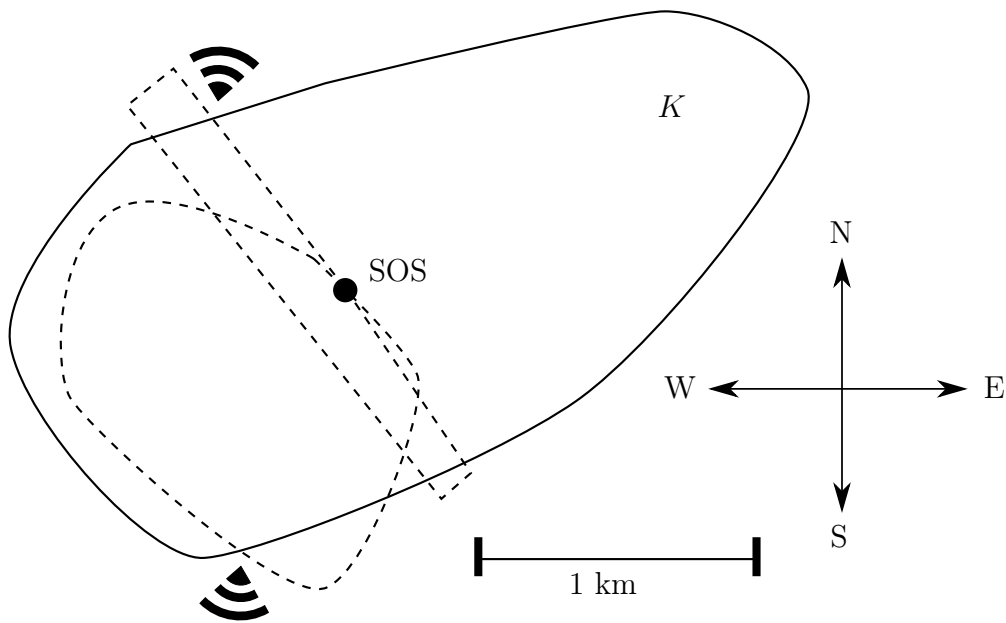
\caption[Visualization of the Minkowski escape problem]{Visualization of the Minkowski escape problem for the special case of two dimensions with Euclidean measurement. This presents two possible Minkowski escape paths which, however, are not the length-minimizing one. For this $K$, the shortest Minkowski escape path is most likely a closed polygonal cruve with two vertices.}
\label{img:minkowskiescapeproblem}
\end{figure}

In order to indicate the connection to Minkowski worm problems in our setting, we will call the problem the \textit{Minkowski escape problem}. We start by generalizing the problem to any dimension. So, we are considering higher dimensional \say{forests} which one aims to escape. We let $K\subset\R^n$ be a convex body, measure lengths by $\ell_T$, where $T\subset\R^n$ is a convex body, and we call $\gamma$ a \textit{closed Minkowski escape path} for $K$ if $\gamma$ is a closed curve and for each point $P\in K$ and each translation $\mu$ for which $P=\mu(\gamma(0))$ the intersection $\mu\left(\{\gamma\}\right)\cap \partial K$ is nonempty. So, in contrast to considering not necessarily closed paths, allowing the motions to be Euclidean motions and measuring the lengths in the standard Euclidean sense in the escape problem of Finch and Wetzel, we only consider closed paths, translations and measure the lengths by the metric induced by the Minkowski functional with respect to the polar of $T$. Translating this problem into \say{our (mesocosmic) reality}--therefore, requiring $n=2$ and Euclidean measurements--, we get a slightly different problem (of course there are no limits to creativity) (cf.\;Figure \ref{img:minkowskiescapeproblem}):\\

\par
\begingroup
\leftskip=1cm
\noindent \textit{Two hikers walk in a forest. One of them gets injured and is in need of medical attention. The unharmed hiker would like to make the emergency call. Although he has his cell phone with him, there is only reception outside the forest. He has a map of the forest, i.e., the shape of the forest and its dimensions are known to him, and a compass to orient himself in terms of direction. Furthermore, he is able to measure the distance he has walked. However, he does not know exactly where in the forest he is. What's the best way to get out of the forest, put off the emergency call, and then get back to the injured hiker?}\\
\par
\endgroup

\noindent The fact that in our story the unharmed hiker knows the shape of the forest and has a compass to orient himself in terms of direction is due to the fact that in our Minkowski escape problem, translations are the only allowed motions. The condition of coming back to the injured hiker is a consequence of our demand to consider only closed curves.

We can prove the analogue to Theorem \ref{Thm:FinchWetzel}:

\bthm\label{Thm:analogueFinchWetzel}
Let $K,T\subset\R^n$ be convex bodies. Then, an/the $\ell_T$-minimizing closed Minkowski escape path for $K$ has $\ell_T$-length $\alpha^*$ if and only if $\alpha^*$ is the largest $\alpha$ for which
\beqq K\in A(T,\alpha),\eeqq
i.e., for which for every closed path $\gamma$ of $\ell_T$-length $\leq \alpha$, there is a translation $\mu$ such that $K$ covers $\mu\left(\{\gamma\}\right)$.
\ethm

Having in mind that Minkowski escape paths for a convex body $K\subset\R^n$ can be understood as closed curves which cannot be translated into the interior of $K$, we can use the Minkowski billiard characterization of shortest closed polygonal curves that cannot be translated into the interior of $K$, in order to directly conclude the following corollary. Note for this line of argumentation that shortest closed curves that cannot be translated into the interior of $K$ are in fact closed polygonal curves.

\bcor\label{Cor:analogue}
Let $K,T\subset\R^n$ be convex bodies, where $T$ is additionally assumed to be strictly convex. An/The $\ell_T$-minimizing closed $(K,T)$-Minkowski billiard trajectory has $\ell_T$-length $\alpha^*$ if and only if $\alpha^*$ is the largest $\alpha$ for which
\beqq K\in A(T,\alpha).\eeqq
\ecor

So, the unharmed hiker in our story can conceptualize his problem by searching for length-minimizing closed Euclidean billiard trajectories.

In general, the problem of minimizing over Minkowski escape problems in the sense of varying the forest while maintaining their volume in order to find the forest with minimal escape length becomes the problem of solving systolic Minkowski billiard inequalities, or equivalently, the problem of proving/investigating Viterbo's conjecture for Lagrangian products in $\R^n\times\R^n$.

This means: If the hikers want to play it safe from the outset by choosing, among forests of equal area, the one where the time needed to help an injured hiker is minimized, then it is useful for them to be familiar with symplectic geometry or billiard dynamics. Of course, they could have paid attention from the beginning to where they entered the forest from and how they designed their path. Then they do not have to solve too difficult problems.

This paper is organized as follows: In Section \ref{Sec:Preliminaries}, we start with some relevant preliminaries before, in Section \ref{Sec:Properties}, we derive properties of Minkowski worm problems and the fundamental results in order to prove Theorems \ref{Thm:main1}, \ref{Thm:Mahler}, \ref{Thm:relations}, \ref{Thm:main2}, and \ref{Thm:main3} and Corollary \ref{Cor:optimizationmain1} in Sections \ref{Sec:main1}, \ref{Sec:main2}, and \ref{Sec:main3}. In Section \ref{Sec:Wetzelproblem}, we prove that it is justified to transfer a special case of Viterbo's conjecture into one for Wetzel's problem which becomes Conjecture \ref{Conj:Wetzelsproblem}. In Section \ref{Sec:escape}, we prove Theorem \ref{Thm:FinchWetzel} as analogue to the relationship between Moser's worm problem and Bellman's lost-in-a-forest problem. Finally, in Section \ref{Sec:OptimizationProblem}, we discuss a computational approach for improving lower bounds in Minkowski worm problems, especially lower bounds for Wetzel's problem.

\section{Preliminaries}\label{Sec:Preliminaries}

We begin by collecting some results concerning the \textit{Fenchel conjugate} of a convex and continuous function $H:\R^n\rightarrow\R$, which for $x^*\in\R^n$ is defined by
\beqq  H^*(x^*)=\sup_{x\in \R^n}(\langle x,x^*\rangle - H(x)).\eeqq

\bprop[Proposition II.1.8 in \cite{Ekeland1990}]\label{Prop:ConsofFenchelConjugateDef}
If $H^*$ is the Fenchel conjugate of a convex and continuous function $H:\R^n\rightarrow\R$, then for $x\in\R^n$ we have
\beqq H(x)=\sup_{x^*\in\R^n}(\langle x^*,x\rangle - H^*(x^*)).\eeqq
\eprop

The \textit{subdifferential} of $H$ in $x\in\R^n$ is given by
\beqq \partial H(x)=\{x^*\in\R^n \vert H^*(x^*)=\langle x^*,x\rangle -H(x)\}\eeqq

Then, we get the \textit{Legendre recipocity formula}:

\bprop[Proposition II.1.15 in \cite{Ekeland1990}]\label{Prop:Legendrereciprocityformula}
For a convex and continuous function $H:\R^n\rightarrow\R$ the Legendre reciprocity formula is given by
\beqq x^*\in\partial H(x)\;\Leftrightarrow\; H^*(x^*)+H(x)=\langle x^*,x\rangle \;\Leftrightarrow\; x\in\partial H^*(x^*),\eeqq
where $x,x^*\in\R^n$.
\eprop

We state the \textit{generalized Euler identity}:

\bprop\label{Prop:Euleridentity}
Let $H:\R^n\rightarrow \R$ be a $p$-positively homogeneous, convex and continuous function of $\R^n$. Then, for each $x\in \R^n$ the following identity holds:
\beqq \langle x^*,x\rangle=pH(x)\quad \forall x^*\in\partial H(x).\eeqq
\eprop

\bpf
For each $x\in \R^n$, since $H$ is convex and continuous, we have
\beqq \partial H(x)\neq \emptyset.\eeqq
For each
\beqq x^*\in\partial H(x)\eeqq
Proposition \ref{Prop:Legendrereciprocityformula} provides
\beq H^*(x^*)+H(x)=\langle x^*,x\rangle\label{eq:Eulerlem1}\eeq
and from Proposition \ref{Prop:ConsofFenchelConjugateDef}, i.e.,
\beqq H(x)=\sup_{x^*\in \R^n}(\langle x^*,x\rangle -H^*(x^*)),\eeqq
we get
\beq H(y)\geq \langle x^*,y \rangle-H^*(x^*)\quad \forall y\in \R^n.\label{eq:Eulerlem2}\eeq
Combining \eqref{eq:Eulerlem1} and \eqref{eq:Eulerlem2} we get
\beq H(y)\geq \langle x^*,y-x\rangle +H(x)\quad \forall y\in \R^n.\label{eq:Eulerlem3}\eeq
Now, we set
\beqq y=\lambda x \quad (\lambda >0)\eeqq
and recognize to have equality in \eqref{eq:Eulerlem3} for $\lambda \rightarrow 1$. Furthermore, we obtain by the $p$-homogeneity of $H$ for $\lambda \rightarrow 1$:
\beqq \lim_{\lambda \rightarrow 1}\frac{g(\lambda)-g(1)}{\lambda -1}H(x)= \langle x^*,x\rangle,\eeqq
where we introduced the function
\beqq g(x):=x^p.\eeqq
Because of
\beqq g'(1)=p\eeqq
we get
\beqq pH(x)=\langle x^*,x\rangle.\eeqq
\epf

Noting that for convex body $C\subset\R^n$
\beqq H_C=\frac{1}{2}\mu_C^2\eeqq
is a $2$-positively homogeneous, convex and continuous function, we derive the following properties:

\bprop\label{Lem:ConsLegendre}
For convex body $C\subset\R^n$ we have
\beqq H_C^*=H_{C^\circ}.\eeqq
\eprop

\bpf
For $\xi\in\R^{n}$ we have
{\allowdisplaybreaks\begin{align*}
\mu_{C^{\circ}}(\xi)&=\min\lbrace t\geq 0 : \xi \in tC^{\circ}\rbrace\\
&=\min \left\lbrace t\geq 0 : \xi \in t \lbrace \widehat{\xi}\in\R^{n}: \langle\widehat{\xi},x\rangle\leq 1 \,\forall x\in C\rbrace\right\rbrace\\
&=\min \left\lbrace t\geq 0 : \xi \in \lbrace \widehat{\xi}\in\R^{n}: \langle \widehat{\xi},x\rangle\leq t \,\forall x\in C\rbrace\right\rbrace\\
&= \min \lbrace t\geq 0 : \langle \xi,x \rangle \leq t \,\forall x\in C\rbrace\\
&=\max_{x\in C}\;\langle \xi,x\rangle\\
& =\max_{\mu_C(x)=1}\langle \xi,x\rangle,
\end{align*}}%
and therefore
{\allowdisplaybreaks\begin{align*}
H_C^*(\xi)&=\sup_{x\in\R^{n}}\left(\langle\xi,x\rangle-H_C(x)\right)\\
&=\sup_{r\geq 0}\sup_{\mu_C(x)=1}\left( \langle\xi,rx\rangle-\frac{1}{2}\mu_C(rx)^2\right)\\
&= \sup_{r\geq 0} \left( r \left(\sup_{\mu_C(x)=1}\langle\xi,x\rangle\right) - \frac{r^2}{2}\right) \\
& = \max_{r\geq 0} \left( r \left(\max_{\mu_C(x)=1}\langle\xi,x\rangle\right) - \frac{r^2}{2}\right)\\
&=\max_{r\geq 0}\left(r\mu_{C^{\circ}}(\xi)-\frac{r^2}{2}\right)\\
&=\frac{\mu_{C^{\circ}}(\xi)^2}{2}\\
&=H_{C^{\circ}}(\xi).
\end{align*}}%
\epf

\bprop\label{Prop:ConsEuler}
Let $C\subset\R^n$ be a convex body. If $x^*\in\partial H_C(x)$ for $x\in \R^n$, then
\beqq H_{C^\circ}(x^*)=H_C(x).\eeqq
\eprop

\bpf
With Proposition \ref{Prop:Euleridentity} and the $2$-homogeneity of $H_{C^\circ}$ we can write
\beqq 2H_{C^\circ}(x^*)=\langle x',x^*\rangle,\eeqq
where
\beqq x'\in\partial H_{C^\circ}(x^*),\eeqq
which together with Propositions \ref{Prop:Legendrereciprocityformula} and \ref{Lem:ConsLegendre} and
\beqq (C^\circ)^\circ =C\eeqq
is equivalent to
\beqq x^*\in\partial H_{C^\circ}^*(x')=\partial H_C(x').\eeqq
Therefore, again using Proposition \ref{Prop:Euleridentity}, we can conclude
\beqq 2H_{C^\circ}(x^*)=\langle x',x^*\rangle=2H_{C}(x').\eeqq
In the following we show that
\beqq H_C(x')=H_C(x).\eeqq
This would prove the claim.

Again, using Propositions \ref{Prop:Legendrereciprocityformula} and \ref{Lem:ConsLegendre}, the fact
\beqq x^*\in\partial H_C(x)\eeqq
is equivalent to
\beqq x\in \partial H_C^*(x^*)=\partial H_{C^\circ}(x^*).\eeqq
All previous informations now can be summarized by the following two equations:
\beqq H_C(x)+H_{C^\circ}(x^*)=\langle x,x^*\rangle,\quad H_{C^\circ}(x^*)+H_C(x')=\langle x',x^*\rangle.\eeqq
The difference yields
\beqq H_C(x)-H_C(x')=\langle x-x',x^*\rangle,\eeqq
which implies
\beqq H_C(x')=H_C(x)-\langle x-x',x^*\rangle= H_C(x)-\langle x,x^*\rangle + \langle x',x^*\rangle.\eeqq
The conditions
\beqq x\in\partial H_{C^\circ}(x^*) \; \text{ and } \; x'\in \partial H_{C^\circ}(x^*)\eeqq
imply, applying Proposition \ref{Prop:Euleridentity},
\beqq -\langle x,x^*\rangle + \langle x',x^*\rangle=-2H_{C^\circ}(x^*)+2H_{C^\circ}(x^*)=0,\eeqq
therefore
\beqq H_C(x')=H_C(x).\eeqq
\epf




The following proposition is the generalization of \cite[Proposition 3.11]{KruppRudolf2022} from closed polygonal curves to closed curves:

\bprop\label{Prop:lengthpropertygeneral}
Let $T\subset\R^n$ be a convex body, $q\in cc(\R^n)$ and $\lambda >0$. Then, we have
\beqq \ell_T(\lambda q) = \ell_{\lambda T}(q) = \lambda \ell_T(q).\eeqq
\eprop

\bpf
From
\beqq \mu_{T^\circ}(\lambda x)=\mu_{(\lambda T)^\circ}(x)=\lambda \mu_{T^\circ}(x),\quad x\in\R^n,\eeqq
(cf.\;\cite[Proposition 2.3(iii)]{KruppRudolf2022}) we conclude
\beqq \ell_T(\lambda q) =\int_0^{\widetilde{T}} \mu_{(T)^\circ}((\dot{\lambda q})(t))\,\mathrm{d}t=\int_0^{\widetilde{T}} \mu_{(\lambda T)^\circ}(\dot{q}(t))\,\mathrm{d}t = \ell_{\lambda T}(q)\eeqq
and
\beqq \ell_T(\lambda q)=\int_0^{\widetilde{T}} \mu_{(T)^\circ}((\dot{\lambda q})(t))\,\mathrm{d}t = \lambda \int_0^{\widetilde{T}} \mu_{(T)^\circ}((\dot{q})(t))\,\mathrm{d}t=\lambda \ell_T(q).\eeqq
\epf

We continue by recalling \cite[Theorem 1.1]{Rudolf2022} which will be useful throughout this paper:

\bthm\label{Thm:relationcapacity}
Let $K,T\subset\R^n$ be convex bodies. Then, we have
\beqq c_{EHZ}(K\times T)=\min_{q\in F^{cp}(K)}\ell_T(q)=\min_{p\in F^{cp}(T)}\ell_K(p)=\min_{q \in M_{n+1}(K,T)} \ell_T(q).\eeqq
\ethm

We note that in \cite[Theorem 1.1]{Rudolf2022} actually appear $F^{cp}_{n+1}(K)$ and $F^{cp}_{n+1}(T)$ instead of $F^{cp}(K)$ and $F^{cp}(T)$, respectively. However, for the purposes within this paper, we only need this more general formulation which is valid since there are no $\ell_T$/$\ell_K$-minimizing closed polygonal curves in $F^{cp}(K)$/$F^{cp}(T)$ with more than $n+1$ vertices  and shorter $\ell_T$/$\ell_K$-length than the $\ell_T$/$\ell_K$-minimizing closed polygonal curves in  $F^{cp}_{n+1}(K)$/$F^{cp}_{n+1}(T)$ (cf.\;the proof of point (i) in the proof of \cite[Theorem 2.2]{Rudolf2022}).

We collect some invariance properties of Viterbo's as well as of Mahler's conjecture:

\bprop\label{Prop:invViterbotranslation}
Viterbo's conjecture is invariant under translations.
\eprop

\bpf
Translations
\beqq t_a:\R^n\rightarrow \R^n, \quad \xi \mapsto \xi + a,\; a\in\R^n,\eeqq
are symplectomorphism because of
\beqq \mathrm{d}t_a(\xi)=\mathbb{1}\eeqq
and therefore
\beqq \mathrm{d}t_a(\xi)^TJ\mathrm{d}t_a(\xi)=J.\eeqq

Finally, we know that Viterbo's conjecture is invariant under symplectomorphisms, since symplectomorphisms preserve the volume as well as the action and therefore the EHZ-capacity.
\epf

\bprop\label{Prop:invViterbo}
Let $C\subset\R^{2n}$ and $K,T\subset\R^n$ be convex bodies. Then
\beqq \vol(C)\geq \frac{c_{EHZ}(C)^n}{n!} \;\Leftrightarrow\; \vol(\lambda C)\geq \frac{c_{EHZ}(\lambda C)^n}{n!}\eeqq
for $\lambda >0$, and
\beqq \vol(K \times  T)\geq \frac{c_{EHZ}( K \times T)^n}{n!} \; \Leftrightarrow \; \vol(\lambda K \times \mu T)\geq \frac{c_{EHZ}(\lambda K \times \mu T)^n}{n!}\eeqq
for $\lambda,\mu >0$. If
\beqq \Phi:\R^n\rightarrow\R^n\eeqq
is an invertible linear transformation, then
\begin{align*}
 \vol(K \times  T) &\geq \frac{c_{EHZ}( K \times T)^n}{n!}\\
\Leftrightarrow \; \vol\left(\Phi(K)\times \left(\Phi^T\right)^{-1}(T)\right) &\geq \frac{c_{EHZ}\left(\Phi(K)\times \left(\Phi^T\right)^{-1}(T)\right)^n}{n!}.
\end{align*}
\eprop

\bpf
We have
\beqq \vol(\lambda C)=\lambda^{2n}\vol(C)\eeqq
and
\beqq c_{EHZ}(\lambda C)=\lambda^2 c_{EHZ}(C)\eeqq
due to the $2$-homogeneity of the action. Further,
\beqq \vol(\lambda K \times \mu T)= \vol(\lambda K)\vol(\mu K)=\lambda^n \mu^n \vol(K)\vol(T)=\lambda^n \mu^n \vol(K\times T)\eeqq
and
\beqq c_{EHZ}(\lambda K \times \mu T) = \min_	{q\in F^{cp}(\lambda K)}\ell_{\mu T}(q)= \lambda\mu \min_{q\in F^{cp}(K)}\ell_T(q)\eeqq
due to Theorem \ref{Thm:relationcapacity} and \cite[Proposition 3.11(ii)$\&$(iv)]{KruppRudolf2022} (cf.\;also Lemma \ref{Lem:Fhomogenity}).

Furthermore,
\beqq \Phi\times \left(\Phi^T\right)^{-1}\eeqq
is a symplectomorphism, i.e.,
\beqq \left(\Phi\times \left(\Phi^T\right)^{-1}\right)^T J \left(\Phi\times \left(\Phi^T\right)^{-1}\right) = J.\eeqq

Indeed, for $a,b\in\R^n$, we calculate
\begin{align*}
\left(\Phi\times \left(\Phi^T\right)^{-1}\right)^T J \left(\Phi\times \left(\Phi^T\right)^{-1}\right)(a,b) = & \left(\Phi\times \left(\Phi^T\right)^{-1}\right)^T J \left(\Phi(a),\left(\Phi^T\right)^{-1}(b)\right)\\
=&\left(\Phi\times \left(\Phi^T\right)^{-1}\right)^T \left(\left(\Phi^T\right)^{-1}(b),-\Phi(a)\right)\\
=&\left(\Phi^T\times \Phi^{-1}\right)\left(\left(\Phi^T\right)^{-1}(b),-\Phi(a)\right)\\
=&\left(\Phi^T\left(\left(\Phi^T\right)^{-1}(b)\right),\Phi^{-1}(-\Phi(a))\right)\\
=&(b,-a)\\
=&J(a,b),
\end{align*}
where we used the facts
\beqq \left(\Phi^T\right)^T=\Phi\; \text{ and }\; \left(\Phi^T\right)^{-1}=\left(\Phi^{-1}\right)^{T}.\eeqq

Finally, we recall that every symplectomorphism preserves the volume as well as the action and therefore the EHZ-capacity.

\epf

\bprop\label{Prop:invMahler}
If $T\subset\R^n$ is a centrally symmetric convex body and
\beqq \Phi:\R^n\rightarrow\R^n\eeqq
an invertible linear transformation, then
\beqq \vol(T)\vol(T^\circ)\geq \frac{4^n}{n!}\;\Leftrightarrow \; \vol(\Phi(T))\vol(\Phi(T)^\circ)\geq \frac{4^n}{n!}.\eeqq
\eprop

\bpf
Because of
\beqq \Phi(T)^\circ=\left(\Phi^T\right)^{-1}\left(T^\circ\right)\eeqq
and the volume preservation of
\beqq \Phi\times \left(\Phi^T\right)^{-1},\eeqq
we have
\begin{align*}
\vol\left(\Phi(T)\right)\vol\left((\Phi(T))^\circ\right) =&\vol\left(\Phi(T)\times \Phi(T)^\circ\right)\\
=&\vol\left(\Phi(T)\times \left(\Phi^T\right)^{-1}(T^\circ)\right)\\
=&\vol\left(\left(\Phi\times \left(\Phi^T\right)^{-1}\right)(T\times T^\circ)\right)\\
=&\vol(T\times T^\circ)\\
=&\vol(T)\vol(T^\circ).
\end{align*}
\epf

\section{Properties of Minkowski worm problems}\label{Sec:Properties}

We begin by concluding some basic properties of the set
\beqq A(T,\alpha), \quad T\in\mathcal{C}(\R^n),\; \alpha > 0.\eeqq
We note that all of the following properties can be easily extended to the case $\alpha \geq 0$. Nevertheless, for the sake of simplicity and in order to avoid trivial case distinctions when it is not possible to divide by $\alpha$, for the following we just treat the case $\alpha >0$.

\bprop\label{Prop:Ahomogenity1}
Let $T\subset\R^n$ be a convex body and $\alpha,\lambda,\mu > 0$. Then we have
\beqq A(\lambda T,\mu \alpha) = \frac{\mu}{\lambda}A(T,\alpha).\eeqq
\eprop

\bpf
We have
\beqq A(\lambda T,\mu\alpha)=\{K\in\mathcal{C}(\R^n) : L_{\lambda T}(\mu\alpha) \subseteq C(K)\}.\eeqq
Because of
\beqq \ell_{\lambda T}(q)=\ell_T(\lambda q)\eeqq
(cf.\;Proposition \ref{Prop:lengthpropertygeneral}) we conclude
\beqq q\in L_{\lambda T}(\mu\alpha) \; \Leftrightarrow \; \lambda q \in L_T(\mu\alpha)\eeqq
which together with
\beqq q\in C(K)\; \Leftrightarrow \; \lambda q \in C(\lambda K)\eeqq
implies
\begin{align*}
A(\lambda T,\mu\alpha)=&\{K\in\mathcal{C}(\R^n):q\in L_{\lambda T}(\mu\alpha)\Rightarrow q\in C(K)\}\\
=&\{K\in\mathcal{C}(\R^n):\lambda q\in L_{T}(\mu\alpha)\Rightarrow \lambda q\in C(\lambda K)\}\\
\overset{(K^*=\lambda K)}{\underset{(q^*=\lambda q)}{=}}&\left\{\frac{1}{\lambda}K^*\in\mathcal{C}(\R^n):q^*\in L_{T}(\mu\alpha)\Rightarrow q^*\in C(K^*)\right\}\\
=&\frac{1}{\lambda}A(T,\mu\alpha).
\end{align*}
Again referring to Proposition \ref{Prop:lengthpropertygeneral} we conclude
\beqq \ell_T(q)=\mu\alpha \; \Leftrightarrow\; \ell_T\left(\frac{q}{\mu}\right)=\alpha,\eeqq
and therefore
\beqq q\in L_T(\mu\alpha)\; \Leftrightarrow \; \frac{q}{\mu} \in L_T(\alpha).\eeqq
This implies
\begin{align*}
A(T,\mu\alpha) =& \{K\in\mathcal{C}(\R^n):q\in L_{T}(\mu\alpha)\Rightarrow q\in C(K)\}\\
=& \left\{K\in\mathcal{C}(\R^n):\frac{q}{\mu}\in L_{T}(\alpha)\Rightarrow \frac{q}{\mu}\in C\left(\frac{K}{\mu}\right)\right\}\\
\overset{(K^*=\frac{K}{\mu})}{\underset{(q^*=\frac{q}{\mu})}{=}}& \{\mu K^* \in\mathcal{C}(\R^n) :q^*\in L_T(\alpha)\Rightarrow q^*\in C(K^*)\}\\
=& \mu A(T,\alpha).
\end{align*}
\epf

\bprop\label{Prop:cases1}
Let $T\subset\R^n$ be a convex body and $\alpha_1,\alpha_2 >0$. Then, we have
\beqq \alpha_1 \begin{rcases}\begin{dcases} \leq \\ < \\ = \end{dcases}\end{rcases} \alpha_2 \; \Rightarrow \; A(T,\alpha_1) \begin{rcases}\begin{dcases} \supseteq \\ \supsetneq \\ =  \end{dcases}\end{rcases} A(T,\alpha_2).\eeqq
\eprop

\bpf
We find $\mu >0$ such that
\beqq \mu\alpha_1 = \alpha_2.\eeqq
Then, using Proposition \ref{Prop:Ahomogenity1} we have
\beq A(T,\alpha_2)=A(T,\mu \alpha_1) = \mu A(T,\alpha_1).\label{eq:cases10}\eeq
This implies
\beqq \alpha_1 \begin{rcases}\begin{dcases} \leq \\ < \\ =\end{dcases}\end{rcases} \alpha_2 \; \Leftrightarrow \; \mu \begin{rcases}\begin{dcases} \geq \\ > \\ =\end{dcases}\end{rcases}1\; \Leftrightarrow \; A(T,\alpha_1) \begin{rcases}\begin{dcases} \supseteq \\ \supsetneq \\ = \end{dcases}\end{rcases} A(T,\alpha_2),\eeqq
where the last equivalence follows from the following considerations: If we have \eqref{eq:cases10} with $\mu \geq  1$, then
\beqq K\in A(T,\alpha_2)=\mu A(T,\alpha_1)\eeqq
means that
\beqq \frac{1}{\mu}K\in A(T,\alpha_1),\eeqq
i.e.,
\beqq L_T(\alpha_1)\subseteq C\left(\frac{1}{\mu}K\right) \subseteq C(K).\eeqq
This implies
\beqq K\in A(T,\alpha_1)\eeqq
and therefore 
\beqq A(T,\alpha_2) \subseteq A(T,\alpha_1).\eeqq
The case $\mu > 1$ follows by considering that in this case there can be find a convex body $K^*$ with
\beqq K^*\in A(T,\alpha_1)\setminus A(T,\alpha_2).\eeqq

Indeed, for
\beqq K\in A(T,\alpha_2)\eeqq
we define
\beqq \widehat{K}:=\widehat{\lambda}K,\quad \widehat{\lambda}:=\min\{0<\lambda \leq 1 : \lambda K\in A(T,\alpha_2)\}.\eeqq
Then, one has
\beqq \widehat{K}\in A(T,\alpha_2)\eeqq
and with \eqref{eq:cases10}
\beqq \frac{1}{\mu}\widehat{K}\in A(T,\alpha_1).\eeqq
With
\beqq K^*:=\frac{1}{\mu}\widehat{K}\eeqq
it follows
\beqq K^*\in A(T,\alpha_1)\setminus A(T,\alpha_2)\eeqq
by the definition of $\widehat{K}$.
\epf

For convex body $T\subset\R^n$ and $\alpha > 0$ we define the set
\beqq A^{\leq}(T,\alpha)=\left\{K\in \mathcal{C}(\R^n) : L_{T}^{\leq}(\alpha) \subseteq C(K)\right\},\eeqq
where
\beqq L_T^{\leq}(\alpha)=\left\{q\in cc(\R^n) : 0<\ell_T(q)\leq \alpha\right\}=\bigcup_{0 < \widetilde{\alpha} \leq \alpha} L_T\left(\widetilde{\alpha}\right).\eeqq
Then, we have the following identity:

\bprop\label{Prop:Aidentity}
Let $T\subset\R^n$ be a convex body and $\alpha > 0$. Then, we have
\beqq A(T,\alpha)=A^{\leq}(T,\alpha).\eeqq
\eprop

\bpf
By definition it is clear that
\beqq A^\leq(T,\alpha) \subseteq A(T,\alpha).\eeqq

Indeed, if
\beqq K\in A^{\leq}(T,\alpha),\eeqq
then this means
\beqq L_T(\widetilde{\alpha})\subseteq C(K),\quad \text{for all }0<\widetilde{\alpha}\leq \alpha.\eeqq
For $\widetilde{\alpha}=\alpha$ it follows
\beqq L_T(\alpha)\subseteq C(K)\eeqq
and therefore
\beqq K\in A(T,\alpha).\eeqq

Let $0<\widetilde{\alpha}\leq\alpha$. Then, it follows from Proposition \ref{Prop:cases1} that
\beqq A(T,\alpha) \subseteq A\left(T,\widetilde{\alpha}\right),\quad \text{ for all } 0< \widetilde{\alpha}\leq \alpha.\eeqq
This implies
\beqq A(T,\alpha) \subseteq \bigcap_{0<\widetilde{\alpha}\leq \alpha}A\left(T,\widetilde{\alpha}\right)=A^{\leq}(T,\alpha).\eeqq
\epf

\bprop\label{Prop:Aproperty11}
Let $\alpha > 0$ and $T_1,T_2\subset\R^n$ be two convex bodies. Then, we have
\beqq T_1 \subseteq T_2 \; \Rightarrow \; A(T_1,\alpha) \subseteq A(T_2,\alpha).\eeqq
\eprop

\bpf
Let
\beqq T_1\subseteq T_2.\eeqq
If
\beqq K\in A(T_1,\alpha),\eeqq
then it follows from Proposition \ref{Prop:Aidentity} that
\beqq L_{T_1}^\leq(\alpha)\subseteq C(K).\eeqq
Because of
\beqq \ell_{T_1}(q)\leq \ell_{T_2}(q) \; \text{ for all }\; q\in cc(\R^n),\eeqq
as consequence of
\beqq \mu_{T_1^\circ}(x)\leq \mu_{T_2^\circ}(x)\quad \forall x\in\R^n,\eeqq
it follows that
\beqq L_{T_2}^\leq(\alpha) \subseteq L_{T_1}^\leq(\alpha)\eeqq
and therefore
\beqq L_{T_2}^\leq(\alpha) \subseteq C(K).\eeqq
With Proposition \ref{Prop:Aidentity} this implies
\beqq K\in A\left(T_2,\alpha\right).\eeqq
Consequently, it follows
\beqq A(T_1,\alpha) \subseteq A(T_2,\alpha).\eeqq
\epf

\blem\label{Lem:easyLemma}
Let $T\subset\R^n$ be a convex body and $\alpha > 0$. Further, let $K_1,K_2\subset\R^n$ be two convex bodies with
\beqq K_1 \subseteq K_2.\eeqq
Then it holds:
\beqq K_1\in A(T,\alpha) \; \Rightarrow \; K_2\in A(T,\alpha).\eeqq
\elem

\bpf
Let
\beqq K_1\in A(T,\alpha),\eeqq
i.e.,
\beqq L_T(\alpha)\subseteq C(K_1).\eeqq
It obviously holds
\beqq K_1 \subseteq K_2 \; \Rightarrow \; C(K_1)\subseteq C(K_2).\eeqq
Therefore
\beqq L_T(\alpha) \subseteq C(K_2),\eeqq
i.e.,
\beqq K_2 \in A(T,\alpha).\eeqq
\epf

For the next Lemma we recall the following: If $(M,d)$ is a metric space and $P(M)$ the set of all nonempty compact subsets of $M$, then $(P(M),d_H)$ is a metric space, where by $d_H$ we denote the \textit{Hausdorff metric} $d_H$ which for nonempty compact subsets $X,Y$ of $(M,d)$ is defined by
\beqq d_H(X,Y)=\max\left\{\sup_{x\in X}\inf_{y\in Y} d(x,y),\sup_{y\in Y} \inf_{x\in X} d(x,y)\right\}.\eeqq
Then, $(cc(\R^n),d_H)$ is a metric subspace of the complete metric space $(P(\R^n),d_H)$ which is induced by the Euclidean space $(\R^n,|\cdot|)$. For convex body $K\subset\R^n$ we consider
\beqq (F^{cc}(K),d_H) \; \text{ and } \; (C(K),d_H)\eeqq
as metric subspaces of $(cc(\R^n),d_H)$. We have that
\beqq F^{cc}(K)\setminus C(K) \; \text{ and } \; C(K)\eeqq
are complements of each other in $cc(\R^n)$.

\blem\label{Lem:C(T)closed}
Let $K\subset\R^n$ be a convex body. Then, $(C(K),d_H)$ is a closed metric subspace of $(cc(\R^n),d_H)$. 
\elem

\bpf
Since $C(K)$ is a subset of $cc(\R^n)$, $(C(K),d_H)$ is a metric subspace of the metric space $(cc(\R^n),d_H)$. It remains to show that $C(K)$ is a closed subset of $cc(\R^n)$. For that let $(q_j)_{j\in\N}$ be a sequence of closed curves in $C(K)$ $d_H$-converging to the closed curve $q^*$. If
\beqq q^*\notin C(K),\eeqq
then $q^*$ cannot be translated into $K$. Using the closedness of $K$ in $\R^n$ this means
\beqq \min_{k\in\R^n} \; d_H\left(\partial \conv\{K+k,q^*\},\partial (K+k)\right)=:m>0.\eeqq
Then, we can find a $j_0\in\N$ such that
\beqq d_H(p_j,q^*)<m\eeqq
for all $j\geq j_0$. But this implies
\beqq \min_{k\in\R^n} \; d_H\left(\partial \conv\{K+k,q_j\},\partial (K+k)\right) >0\eeqq
for all $j\geq j_0$, i.e., $p_j$ cannot be translated into $K$ for all $j\geq j_0$, a contradiction to
\beqq q_j\in C(K)\quad \forall j\in\N.\eeqq
Therefore, it follows
\beqq q^*\in C(K),\eeqq
and consequently, $(C(K),d_H)$ is a closed metric subspace of $(cc(\R^n),d_H)$.
\epf

\blem\label{Lem:easyLemma2}
Let $T\subset\R^n$ be a convex body and $(\alpha_k)_{k\in\N}$ an increasing sequence of positive real numbers converging to $\alpha > 0$ for $k\rightarrow\infty$. If
\beqq K\in A(T,\alpha_k)\quad \forall k\in\N,\eeqq
then also
\beqq K\in A(T,\alpha).\eeqq
\elem

\bpf
Let
\beqq K\in A(T,\alpha_k)\quad \forall k\in\N,\eeqq
i.e.,
\beq L_T(\alpha_k)\subseteq C(K)\quad \forall k\in\N.\label{eq:easyLemma1}\eeq
This means for all $k\in\N$ that for all
\beqq q\in cc(\R^n) \; \text{ with } \; \ell_T(q)=\alpha_k\eeqq
holds
\beqq q\in C(K).\eeqq

Let us assume that
\beqq K\notin A(T,\alpha),\eeqq
i.e.,
\beqq L_T(\alpha) \nsubseteq C(K).\eeqq
This means that there is a $q^*\in cc(\R^n)$ with
\beqq \ell_T(q^*)=\alpha \; \text{ and } \; q^*\in F^{cc}(K)\setminus C(K).\eeqq
Due to Lemma \ref{Lem:C(T)closed} $(C(K),d_H)$ is a closed metric subspace of $(cc(\R^n),d_H)$. Since $F^{cc}(K)\setminus C(K)$ is the complement of $C(K)$ in $cc(\R^n)$ it follows that $F^{cc}(K)\setminus C(K)$ is an open subset of the metric space $(cc(\R^n),d_H)$. Consequently there is a $k_0\in \N$ sufficiently big such that
\beq q^*\frac{\alpha_k}{\alpha}\in F^{cc}(K)\setminus C(K)\quad \forall k\geq k_0.\label{eq:easyLemma2}\eeq
But with \cite[Proposition 3.11(iv)]{KruppRudolf2022} it is
\beqq \ell_T\left(q^*\frac{\alpha_k}{\alpha}\right)=\ell_T\left(q^*\right)\frac{\alpha_k}{\alpha}=\alpha_k\quad \forall k\in\N,\eeqq
i.e.,
\beqq q^*\frac{\alpha_k}{\alpha} \in L_T(\alpha_k)\quad \forall k\in\N,\eeqq
which produces a contradiction between \eqref{eq:easyLemma1} and \eqref{eq:easyLemma2}.

Therefore it follows
\beqq K\in A(T,\alpha).\eeqq
\epf

The next Proposition justifies to write \say{$\min$} instead of \say{$\inf$} within the Minkowski worm problem. The main ingredient of its proof will be Blaschke's selection theorem (cf.\;\cite[\S 18]{Blaschke1916}).

\bthm[Blaschke selection theorem]\label{Thm:Blaschkeselection}
Let $(C_k)_{k\in\N}$ be a sequence of convex bodies in $\R^n$ satisfying
\beqq C_k \subset B_R^n,\quad  R>0,\eeqq
for all $k\in\N$. Then there is a subsequence $(C_{k_l})_{l\in\N}$ and a convex body $C$ in $\R^n$ such that $C_{k_l}$ $d_H$-converges to $C$ for $l\rightarrow\infty$. 
\ethm

\bprop\label{Prop:minimumexists}
Let $T$ be a convex body and $\alpha > 0$. Then we have
\beqq \inf_{K\in A(T,\alpha)} \vol(K)=\min_{K\in A(T,\alpha)}\vol(K).\eeqq
\eprop

\bpf
Let $(K_k)_{k\in\N}$ be a minimizing sequence of
\beq \inf_{K\in A(T,\alpha)}\vol(K).\label{eq:minimumexists0}\eeq
Then, there is a $k_0\in\N$ and a sufficiently big $R>0$ such that
\beqq K_k \subset B_R^n \quad \forall k\geq k_0.\eeqq

Indeed, if this is not the case, then there is a subsequence $\left(K_{k_j}\right)_{j\in\N}$ such that
\beq R_j:=\max\{R>0 : K_{k_j} \in F(B_R^n)\} \rightarrow \infty\quad (j\rightarrow\infty).\label{eq:minimumexists1}\eeq
Guaranteeing
\beqq K_{k_j}\in A(T,\alpha)=\{K\in\mathcal{C}(\R^n):L_T(\alpha)\subseteq C(K)\}\quad \forall j\in\N\eeqq
means that
\beq V_j:= \vol\left(K_{k_j}\right) \rightarrow \infty \quad (j\rightarrow\infty).\label{eq:minimumexists2}\eeq
The latter follows together with \eqref{eq:minimumexists1} and the convexity of $K_{k_j}$ for all $j\in\N$ from the fact that due to
\beqq L_T(\alpha)\subseteq C\left(K_{k_j}\right)\quad \forall j\in\N\eeqq
there is no direction in which $K_{k_j}$ can be shrunk. But \eqref{eq:minimumexists2} is not possible since $\left(K_{k_j}\right)_{j\in\N}$ is a minimizing sequence of \eqref{eq:minimumexists0}.

Applying Theorem \ref{Thm:Blaschkeselection}, there is a subsequence $(K_{k_l})_{l\in\N}$ and a convex body $K\subset\R^n$ such that $K_{k_l}$ $d_H$-converges to $K$ for $l\rightarrow\infty$. It remains to show that
\beqq K\in A(T,\alpha).\eeqq

The fact that $K_{k_l}$ $d_H$-converges to $K$ for $l\rightarrow \infty$ implies that for every $\eps >0$ there is an $l_0\in\N$ such that
\beq (1-\eps)K \subseteq K_{k_l} \subseteq (1+\eps)K\quad \forall l\geq l_0.\label{eq:minimumexists3}\eeq
Then, with
\beqq K_{k_l}\in A(T,\alpha)\quad \forall l\in\N\eeqq
it follows from the second inclusion in \eqref{eq:minimumexists3} together with Lemma \ref{Lem:easyLemma} that
\beqq (1+\eps)K \in A(T,\alpha).\eeqq
Applying Proposition \ref{Prop:Ahomogenity1} this means
\beqq K\in A\left(T,\frac{\alpha}{1+\eps}\right).\eeqq

We define the sequence
\beqq \alpha_k:=\frac{\alpha}{1+\frac{1}{k}} \quad \forall k\in\N.\eeqq
Then, $(\alpha_k)_{k\in\N}$ is an increasing sequence of positive numbers converging to $\alpha$ for $k\rightarrow \infty$ and together with the aboved mentioned ($\eps >0$ can be chosen arbitrarily) we have
\beqq K\in A\left(T,\alpha_k\right)\quad \forall k\in\N.\eeqq
Applying Lemma \ref{Lem:easyLemma2} it follows
\beqq K\in A(T,\alpha).\eeqq
\epf

\bprop\label{Prop:Ahomogenity2}
Let $T\subset\R^n$ be a convex body and $\alpha,\lambda,\mu > 0$. Then we have
\beqq \min_{K\in A(\lambda T,\mu \alpha)}\vol(K)=\frac{\mu^n}{\lambda^n} \min_{K\in A(T,\alpha)}\vol(K).\eeqq
\eprop

\bpf
From
\beqq A(\lambda T,\mu\alpha)=\frac{1}{\lambda} A(T,\mu \alpha)\eeqq
(cf.\;Proposition \ref{Prop:Ahomogenity1}) it follows
\beqq K\in A(\lambda T,\mu\alpha) \;\Leftrightarrow\; \lambda K\in A(T,\mu \alpha)\eeqq
and therefore
\beqq \min_{K\in A(\lambda T,\mu\alpha)}\vol(K)\overset{(K^*=\lambda K)}{=}\min_{K^*\in A(T,\mu\alpha)}\vol\left(\frac{K^*}{\lambda}\right)=\frac{1}{\lambda^n} \min_{K^*\in A(T,\mu\alpha)}\vol(K^*).\eeqq
From
\beqq A(T,\mu\alpha)=\mu A(T,\alpha)\eeqq
(cf.\;Proposition \ref{Prop:Ahomogenity1}) it follows
\beqq K\in A(T,\mu\alpha) \;\Leftrightarrow\; \frac{K}{\mu}\in A(T,\alpha)\eeqq
and therefore
\beqq \min_{K\in A(T,\mu\alpha)}\vol(K)\overset{(K^*=\frac{K}{\mu})}{=}\min_{K^*\in A(T,\alpha)}\vol(\mu K^*)=\mu^n \min_{K^*\in A(T,\alpha)}\vol(K^*).\eeqq
\epf

\bprop\label{Prop:monotony1}
Let $T\subset\R^n$ be a convex body and $\alpha_1,\alpha_2 >0$. Then, we have
\beqq \alpha_1 \begin{rcases}\begin{dcases} \leq \\ < \\ = \end{dcases}\end{rcases} \alpha_2 \; \Leftrightarrow \; \min_{K\in A(T,\alpha_1)}\vol(K) \begin{rcases}\begin{dcases} \leq \\ < \\ =  \end{dcases}\end{rcases} \min_{K\in A(T,\alpha_2)}\vol(K)\eeqq
\eprop

\bpf
We find $\mu >0$ such that
\beqq \mu\alpha_1 = \alpha_2.\eeqq
Then, we apply Proposition \ref{Prop:Ahomogenity2}.
\epf

Now, for convex bodies $K,T\subset\R^n$ we will turn our attention to the minimization problem
\beqq \min_{q\in F^{cp}(K)}\ell_T(q).\eeqq
The existence of the minimum is guaranteed by Theorem \ref{Thm:relationcapacity}.

\blem\label{Lem:Fhomogenity}
Let $K,T\subset\R^n$ be convex bodies and $\lambda >0$. Then
\beqq \min_{q\in F^{cp}(\lambda K)}\ell_T(q)=\min_{q\in F^{cp}(K)}\ell_T(\lambda q)=\lambda\min_{q\in F^{cp}(K)}\ell_T(q).\eeqq
\elem

\bpf
Similar to \cite[Proposition 3.11(ii)]{KruppRudolf2022} we have
\beqq q\in F^{cp}(\lambda K) \; \Leftrightarrow \; \frac{q}{\lambda}\in F^{cp}(K)\eeqq
and using \cite[Proposition 3.11(iv)]{KruppRudolf2022} therefore
\beqq \min_{q\in F^{cp}(\lambda K)}\ell_T(q)=\min_{\frac{q}{\lambda}\in F^{cp}(K)}\ell_T(q)\overset{(q^*=\frac{q}{\lambda})}{=}\min_{q^*\in F^{cp}(K)}\ell_T(\lambda q^*)=\lambda \min_{q^*\in F^{cp}(K)}\ell_T(q^*).\eeqq
\epf

In the following for convex body $T\subset\R^n$ and $c>0$ we consider the minimax problem\footnote{Whenever we write \beqq \max_{\vol(K)=c}\;\min_{q\in F^{cp}(K)}\ell_T(q)\eeqq the maximum is understood to consider only convex bodies $K\subset\R^n$. This is implicitly indicated by the fact that we defined $F^{cp}(K)$ only for convex bodies $K\subset\R^n$.}
\beqq \max_{\vol(K)=c}\;\min_{q\in F^{cp}(K)}\ell_T(q).\eeqq
The following Proposition guarantees the existence of its maximum:

\bprop\label{Prop:minimaxproblem1}
Let $T\subset\R^n$ be a convex body and $c>0$. Then, we have
\beqq \sup_{\vol(K)=c}\;\min_{q\in F^{cp}(K)}\ell_T(q)=\max_{\vol(K)=c}\;\min_{q\in F^{cp}(K)}\ell_T(q).\eeqq
\eprop

\bpf
Let $(K_k)_{k\in\N}$ be a maximizing sequence of
\beq \sup_{\vol(K)=c}\;\min_{q\in F^{cp}(K)}\ell_T(q).\label{eq:minimaxproblem10}\eeq
Then, there is a $k_0\in\N$ and an $R>0$ such that
\beq K_k \subset B_R^n\quad \forall k\geq k_0.\label{eq:minimaxproblem1}\eeq

Indeed, if this is not the case, then there is a subsequence $\left(K_{k_j}\right)_{j\in\N}$ such that
\beq R_j:=\max\{R>0 : K_{k_j} \in F(B_R^n)\}\rightarrow \infty \quad (j\rightarrow \infty).\label{eq:minimaxproblem11}\eeq
But this implies
\beq L_j:=\min\left\{\ell_T(q) : q\in F^{cp}(K_{k_j})\right\}\rightarrow 0 \quad (j\rightarrow\infty).\label{eq:minimaxproblem12}\eeq
This follows from the fact that for every $j\in\N$ we can find a
\beqq q_j\in F^{cp}(K_{k_j})\eeqq
with
\beqq \ell_T(q_j)\rightarrow 0 \quad (j\rightarrow\infty).\eeqq
The latter is a consequence of \eqref{eq:minimaxproblem11} and the constraint
\beq \vol\left(K_{k_j}\right)=c\quad \forall j\in\N,\label{eq:minimaxproblem13}\eeq
i.e., due to the convexity of $K_{k_j}$ for all $j\in\N$ guaranteeing \eqref{eq:minimaxproblem13} there are directions from the origin in which $K_{k_j}$ has to shrink for $j\rightarrow\infty$ and which are suitable in order to construct convenient $q_j$. But \eqref{eq:minimaxproblem12} is not possible since $\left(K_{k_j}\right)_{j\in\N}$ is a maximizing sequence of \eqref{eq:minimaxproblem10}.

Then, we apply Theorem \ref{Thm:Blaschkeselection} and find a subsequence $(K_{k_l})_{l\in\N}$ and a convex body $K\subset\R^n$ such that $K_{k_l}$ $d_H$-converges to $K$ for $l\rightarrow\infty$. It remains to prove that
\beqq \vol(K)=c,\eeqq
but this is an immediate consequence of the $d_H$-continuity of the volume function.
\epf

\bprop\label{Prop:monotony2}
Let $T\subset\R^n$ be a convex body. Then,
\beqq \max_{\vol(K)=c}\;\min_{q\in F^{cp}(K)}\ell_T(q)\eeqq
increases/decreases strictly if and only if this is the case for $c> 0$.
\eprop

\bpf
We make use of the implication
\beq \max_{\vol(K)=c_1}\;\min_{q\in F^{cp}(K)}\ell_T(q)=\max_{\vol(K)=c_2}\;\min_{q\in F^{cp}(K)}\ell_T(q)\; \Rightarrow\; c_1=c_2\label{eq:monotony21}\eeq
for all $c_1,c_2>0$.

Let us verify \eqref{eq:monotony21}: We assume
\beq \max_{\vol(K)=c_1}\;\min_{q\in F^{cp}(K)}\ell_T(q)=\max_{\vol(K)=c_2}\;\min_{q\in F^{cp}(K)}\ell_T(q)\label{eq:monotony22}\eeq
and without loss of generality $c_1 < c_2$. Let the pair
\beqq (K_1^*,q_1^*) \;  \text{ with } \; \vol(K_1^*)=c_1 \; \text{ and } \; q_1^*\in F^{cp}(K_1^*)\eeqq
be a maximizer of the left side in \eqref{eq:monotony22}, i.e.,
\beqq \max_{\vol(K)=c_1}\;\min_{q\in F^{cp}(K)}\ell_T(q)=\min_{q\in F^{cp}(K_1^*)}\ell_T(q)=\ell_T(q_1^*).\eeqq
With
\beqq q_1^*\in F^{cp}(K_1^*)\eeqq
similar to \cite[Proposition 3.11(ii)]{KruppRudolf2022} we have
\beqq \sqrt[n]{\frac{c_2}{c_1}}q_1^*\in F^{cp}\left(\widetilde{K}\right)\eeqq
for
\beqq \widetilde{K}:=\sqrt[n]{\frac{c_2}{c_1}}K_1^*.\eeqq
From
\beqq \min_{q\in F^{cp}(K_1^*)}\ell_T(q)=\ell_T(q_1^*)\eeqq
it follows together with Lemma \ref{Lem:Fhomogenity} that
\beqq \min_{q\in F^{cp}\left(\widetilde{K}\right)}\ell_T(q)=\min_{q\in F^{cp}\left(\sqrt[n]{\frac{c_2}{c_1}}K_1^*\right)}\ell_T(q)=\sqrt[n]{\frac{c_2}{c_1}} \min_{q\in F^{cp}(K_1^*)}\ell_T(q)=\sqrt[n]{\frac{c_2}{c_1}} \ell_T(q_1^*).\eeqq
Since
\beqq \vol\left(\widetilde{K}\right)=\vol\left(\sqrt[n]{\frac{c_2}{c_1}}K_1^*\right)=\frac{c_2}{c_1}\vol(K_1^*)=c_2,\eeqq
we conclude
\begin{align*}
\max_{\vol(K)=c_1}\;\min_{q\in F^{cp}(K)}\ell_T(q)=\min_{q\in F^{cp}(K_1^*)}\ell_T(q)&=\ell_T(q_1^*)\\
&<\sqrt[n]{\frac{c_2}{c_1}}\ell_T(q_1^*)\\
&=\min_{q\in F^{cp}\left(\widetilde{K}\right)}\ell_T(q)\\
&\leq \max_{\vol(K)=c_2}\;\min_{q\in F^{cp}(K)}\ell_T(q),
\end{align*}
which is a contradiction to \eqref{eq:monotony22}. Therefore, noting that the assumption $c_1 > c_2$ would have led analogously to the same contradiction, it follows
\beqq c_1=c_2.\eeqq

We now prove the equivalence
\beq \max_{\vol(K)=c_1}\;\min_{q\in F^{cp}(K)}\ell_T(q) <\max_{\vol(K)=c_2}\;\min_{q\in F^{cp}(K)}\ell_T(q) \; \Leftrightarrow \; c_1 <c_2\label{eq:monotony25}\eeq
for $c_1,c_2>0$.

If
\beq \max_{\vol(K)=c_1}\;\min_{q\in F^{cp}(K)}\ell_T(q) <\max_{\vol(K)=c_2}\;\min_{q\in F^{cp}(K)}\ell_T(q),\label{eq:monotony23}\eeq
then from the first part of the proof it necessarily follows $c_1 \neq c_2$. Let us assume $c_1>c_2$. We further assume that the pair
\beqq (K_2^*,q_2^*) \; \text{ with } \; \vol(K_2^*)=c_2 \; \text{ and } \; q_2^*\in F^{cp}(K_2^*)\eeqq
is a maximizer of the right side in \eqref{eq:monotony23}, i.e.,
\beqq \max_{\vol(K)=c_2}\;\min_{q\in F^{cp}(K)}\ell_T(q)=\min_{q\in K_2^*}\ell_T(q)=\ell_T(q_2^*).\eeqq
We define
\beqq \widehat{K}:=\sqrt[n]{\frac{c_1}{c_2}}K_2^*.\eeqq
From
\beqq \min_{q\in K_2^*}\ell_T(q)=\ell_T(q_2^*)\eeqq
it follows together with Lemma \ref{Lem:Fhomogenity} that
\beqq \min_{q\in F^{cp}\left(\widehat{K}\right)}\ell_T(q)=\min_{q\in F^{cp}\left(\sqrt[n]{\frac{c_1}{c_2}}K_2^*\right)}\ell_T(q)=\sqrt[n]{\frac{c_1}{c_2}} \min_{q\in F^{cp}(K_2^*)}\ell_T(q)=\sqrt[n]{\frac{c_1}{c_2}} \ell_T(q_2^*).\eeqq
Since
\beqq \vol\left(\widehat{K}\right)=\vol\left(\sqrt[n]{\frac{c_1}{c_2}}K_2^*\right)=\frac{c_1}{c_2}\vol(K_2^*)=c_1,\eeqq
we conclude
\begin{align*}
\max_{\vol(K)=c_2}\;\min_{q\in F^{cp}(K)}\ell_T(q)=\min_{q\in F^{cp}(K_2^*)}\ell_T(q)&=\ell_T(q_2^*)\\
&<\sqrt[n]{\frac{c_1}{c_2}}\ell_T(q_2^*)\\
&=\min_{q\in F^{cp}\left(\widehat{K}\right)}\ell_T(q)\\
&\leq \max_{\vol(K)=c_1}\;\min_{q\in F^{cp}(K)}\ell_T(q),
\end{align*}
which is a contradiction to \eqref{eq:monotony23}. Therefore, we conclude $c_1<c_2$.

Conversely, let $c_1<c_2$. From \eqref{eq:monotony21} we conclude
\beq \max_{\vol(K)=c_1}\;\min_{q\in F^{cp}(K)}\ell_T(q) \neq \max_{\vol(K)=c_2}\;\min_{q\in F^{cp}(K)}\ell_T(q).\label{eq:monotony24}\eeq
If the strict inequality "$>$" holds in \eqref{eq:monotony24}, then we conclude from the above proven implication "$\Rightarrow$" in \eqref{eq:monotony25} that $c_1 > c_2$, a contradiction. Therefore, it follows
\beqq \max_{\vol(K)=c_1}\;\min_{q\in F^{cp}(K)}\ell_T(q) < \max_{\vol(K)=c_2}\;\min_{q\in F^{cp}(K)}\ell_T(q).\eeqq
\epf

\bprop\label{Prop:p^*minF}
Let $K,T\subset\R^n$ be convex bodies with $q^*$ as minimizer of
\beqq \min_{q\in F^{cp}(K)}\ell_T(q).\eeqq
Then, it follows
\beqq K\in A(T,\ell_T(q^*))=A\left(T,\min_{q\in F^{cp}(K)}\ell_T(q)\right).\eeqq
\eprop

\bpf
Let $q^*$ be a minimizer of
\beqq \min_{q\in F^{cp}(K)}\ell_T(q).\eeqq
Then, it follows
\beqq L_T(\ell_T(q^*)) \subseteq C(K).\eeqq

Indeed, otherwise, if there is
\beqq \widetilde{q}\in L_T(\ell_T(q^*))\setminus C(K),\eeqq
i.e.,
\beqq \ell_T(\widetilde{q})=\ell_T(q^*) \; \text{ and } \; \widetilde{q}\in F^{cc}(K)\setminus C(K),\eeqq
then, due to the openess of
\beqq F^{cc}(K)\setminus C(K) \; \text{ in } \; cc(\R^n)\eeqq
with respect to $d_H$ (cf.\;Lemma \ref{Lem:C(T)closed}), there is a $\lambda <1$ such that
\beqq \lambda \widetilde{q} \in F^{cc}(K)\setminus C(K).\eeqq
Then, using \cite[Proposition 3.11(iv)]{KruppRudolf2022}, we conclude
\beqq \ell_T(\lambda \widetilde{q})=\lambda \ell_T(\widetilde{q})< \ell_T(\widetilde{q})=\ell_T(q^*)=\min_{q\in F^{cp}(K)}\ell_T(q).\eeqq
Because of the $d_H$-density of $F^{cp}(K)$ in $F^{cc}(K)$ and the $d_H$-continuity of $\ell_T$ on $F^{cc}(K)$ (cf.\;\cite[Proposition 3.11(v)]{KruppRudolf2022}--which is also valid for closed curves) then we can find a
\beqq \widehat{q}\in F^{cp}(K)\eeqq
with
\beqq \ell_T\left(\widehat{q}\right)< \min_{q\in F^{cp}(K)}\ell_T(q),\eeqq
a contradiction.

Finally, from
\beqq L_T(\ell_T(q^*)) \subseteq C(K)\eeqq
it follows
\beqq K\in A(T,\ell_T(q^*))=A\left(T,\min_{q\in F^{cp}(K)}\ell_T(q)\right).\eeqq
\epf

\blem\label{Lem:intersectionFC}
Let $K\subset\R^n$ be a convex body and $\lambda >1$. If
\beq q\in F^{cc}(K)\cap C(K),\label{eq:intersectionFC1}\eeq
then it follows that
\beq \lambda q \in F^{cc}(K)\setminus C(K).\label{eq:intersectionFC2}\eeq
\elem

\bpf
If we assume \eqref{eq:intersectionFC1} but \eqref{eq:intersectionFC2} does not hold. Then it follows
\beqq q,\lambda q \in C(K)\eeqq
and due to $\lambda >1$ therefore
\beqq q\in C\left(\mathring{K}\right).\eeqq
But this is a contradiction to
\beqq q\in F^{cc}(K).\eeqq
Therefore, it follows \eqref{eq:intersectionFC2}.
\epf

\bprop\label{Prop:later}
Let $T\subset\R^n$ be a convex body and $\alpha > 0$. If $K^*$ is a minimizer of
\beq \min_{K\in A(T,\alpha)}\vol(K),\label{eq:later}\eeq
then
\beqq \min_{q\in F^{cp}(K^*)}\ell_T(q)=\alpha.\eeqq
\eprop

\bpf
If $q^*$ is a minimizer of
\beqq \min_{q\in F^{cp}(K^*)}\ell_T(q),\eeqq
then it follows from Proposition \ref{Prop:p^*minF} that
\beqq K^*\in A(T,\ell_T(q^*)).\eeqq
This means
\beqq \min_{K\in A(T,\ell_T(q^*))}\vol(K)\leq \vol(K^*)=\min_{K\in A(T,\alpha)}\vol(K).\eeqq
Proposition \ref{Prop:monotony1} implies
\beqq \ell_T(q^*)\leq \alpha.\eeqq
If
\beqq \ell_T(q^*)<\alpha,\eeqq
then with Proposition \cite[Proposition 3.11(iv)]{KruppRudolf2022} there is $\lambda >1$ such that
\beqq \ell_T(\lambda q^*)=\alpha.\eeqq
Together with
\beqq F^{cp}(K^*)\subseteq F^{cc}(K^*)\eeqq
and Lemma \ref{Lem:intersectionFC} the fact
\beqq q^*\in F^{cp}(K^*)\eeqq
implies
\beqq \lambda q^*\in F^{cc}(K^*)\setminus C(K^*),\eeqq
therefore, there is no translate of $K^*$ that covers $\lambda q^*$. Consequently,
\beqq K^*\notin A(T,\ell_T(\lambda q^*))=A(T,\alpha),\eeqq
a contradiction to the fact that $K^*$ is a minimizer of \eqref{eq:later}. Therefore, it follows that
\beqq \min_{q\in F^{cp}(K^*)}\ell_T(q)=\ell_T(q^*)=\alpha.\eeqq
\epf

The idea which underlies the following Theorem leads to the heart of this paper.

\bthm\label{Thm:MainProperty1}
Let $T\subset\R^n$ be a convex body. If $K^*$ is a minimizer of
\beq \min_{K\in A(T,\alpha)} \vol(K)\label{eq:MoserRelation1}\eeq
for $\alpha > 0$, then $K^*$ is a maximizer of
\beq \max_{\vol(K)=c}\; \min_{q\in F^{cp}(K)} \ell_T(q)\label{eq:MoserRelation2}\eeq
for
\beqq c:=\vol(K^*)\eeqq
with
\beqq \min_{q\in F^{cp}(K^*)}\ell_T(q)=\alpha.\eeqq

Conversely, if $K^*$ is a maximizer of \eqref{eq:MoserRelation2} for $c>0$, then $K^*$ is a minimizer of \eqref{eq:MoserRelation1} for
\beqq \alpha:= \min_{q\in F^{cp}(K^*)}\ell_T(q)\eeqq
and with
\beqq \vol(K^*)=c.\eeqq

Consequently, for $\alpha, c>0$ we have the equivalence
\beqq \min_{K\in A(T,\alpha)}\vol(K)= c \; \Leftrightarrow \; \max_{\vol(K)=c}\;\min_{q\in F^{cp}(K)}\ell_T(q)= \alpha\eeqq
and moreover
\beq \min_{K\in A(T,\alpha)}\vol(K)\geq c \; \Leftrightarrow \; \max_{\vol(K)=c}\;\min_{q\in F^{cp}(K)}\ell_T(q)\leq \alpha.\label{eq:mainproperty11}\eeq
\ethm

\bpf
Let $K^*$ be a minimizer of \eqref{eq:MoserRelation1} for $\alpha>0$. If $K^*$ is not a maximizer of \eqref{eq:MoserRelation2} for
\beqq c=\vol(K^*),\eeqq
then there is a convex body
\beqq K^{**}\subset\R^n \; \text{ with } \; \vol(K^{**})=c\eeqq
and a
\beqq q^{**}\in F^{cp}(K^{**})\eeqq
such that
\beq \ell_T(q^{**})=\min_{q\in F^{cp}(K^{**})}\ell_T(q)>\min_{q\in F^{cp}(K^*)}\ell_T(q)=\ell_T(q^*),\label{eq:mainproperty12}\eeq
where by $q^*$ we denote a minimizer of
\beqq \min_{q\in F^{cp}(K^*)}\ell_T(q).\eeqq
From Proposition \ref{Prop:p^*minF} it follows
\beq K^{**}\in A(T,\ell_T(q^{**})),\label{eq:mainproperty13}\eeq
and further from Proposotion \ref{Prop:later} that
\beq \min_{q\in F^{cp}(K^*)}\ell_T(q)=\ell_T(q^*)=\alpha.\label{eq:mainproperty14}\eeq
From \eqref{eq:mainproperty12}, \eqref{eq:mainproperty13} and \eqref{eq:mainproperty14} together with Proposition \ref{Prop:monotony1} we conclude
\begin{align*}
c=\vol(K^{**})\stackrel{\eqref{eq:mainproperty13}}{\geq} \min_{K\in A(T,\ell_T(q^{**}))}\vol(K)&\stackrel{\eqref{eq:mainproperty12}}{>}\min_{K\in A(T,\ell_T(q^*))}\vol(K)\\
&\stackrel{\eqref{eq:mainproperty14}}{=}\min_{K\in A(T,\alpha)}\vol(K)\\
&=\vol(K^*)\\
&=c,
\end{align*}
a contradiction. Therefore, $K^*$ is a maximizer of \eqref{eq:MoserRelation2} for
\beqq c=\vol(K^*).\eeqq

Conversely, let $K^*$ be a maximizer of \eqref{eq:MoserRelation2} for $c>0$ with
\beqq q^*\in F^{cp}(K^*)\eeqq
such that
\beqq \max_{\vol(K)=c}\;\min_{q\in F^{cp}(K)}\ell_T(q)=\min_{q\in F^{cp}(K^*)}\ell_T(q)=\ell_T(q^*)=:\alpha.\eeqq
Then, from Proposition \ref{Prop:p^*minF} it follows that
\beqq K^*\in A(T,\alpha),\eeqq
and consequently
\beqq c=\vol(K^*)\geq \min_{K\in A(T,\alpha)} \vol(K).\eeqq
If $K^*$ is not a minimizer of \eqref{eq:MoserRelation1} for
\beqq \alpha=\ell_T(q^*),\eeqq
then there is a
\beqq K^{**}\in A(T,\alpha)\eeqq
with
\beq c=\vol(K^*)>\min_{K\in A(T,\alpha)}\vol(K)=\vol(K^{**}).\label{eq:mainproperty15}\eeq
Then, from Proposition \ref{Prop:later} it follows that
\beqq \min_{q\in F^{cp}(K^{**})}\ell_T(q)=\alpha.\eeqq
This implies
\begin{align*}
\max_{\vol(K)=c}\;\min_{q\in F^{cp}(K)}\ell_T(q)= \min_{q\in F^{cp}(K^*)} \ell_T(q)&=\ell_T(q^*)\\
&=\alpha\\
&= \min_{q\in F^{cp}(K^{**})}\ell_T(q)\\
&\leq \max_{\vol(K)=\vol(K^{**})}\;\min_{q\in F^{cp}(K)}\ell_T(q),
\end{align*}
which because of \eqref{eq:mainproperty15} is a contradiction to Proposition \ref{Prop:monotony2}. We conclude that $K^*$ is a minimizer of \eqref{eq:MoserRelation1} for
\beqq \alpha=\ell_T(q^*).\eeqq

From the before proven it clearly follows the equivalence
\beqq \min_{K\in A(T,\alpha)}\vol(K)= c \; \Leftrightarrow \; \max_{\vol(K)=c}\;\min_{q\in F^{cp}(K)}\ell_T(q)= \alpha\eeqq
for $\alpha,c >0$. In order to prove \eqref{eq:mainproperty11} it remains to show
\beqq \min_{K\in A(T,\alpha)}\vol(K)> c \; \Leftrightarrow \; \max_{\vol(K)=c}\;\min_{q\in F^{cp}(K)}\ell_T(q) < \alpha.\eeqq

Let $K^*$ be a minimizer of
\beqq \min_{K\in A(T,\alpha)}\vol(K),\eeqq
where $c>0$ is chosen such that
\beq \vol(K^*)> c.\label{eq:mainproperty16}\eeq
Then we know from the above reasoning that $K^*$ is a maximizer of
\beqq \max_{\vol(K)=\vol(K^*)}\;\min_{q\in F^{cp}(K)}\ell_T(q)\eeqq
with
\beqq \min_{q\in F^{cp}(K^*)}\ell_T(q)=\alpha.\eeqq
From \eqref{eq:mainproperty16} and Proposition \ref{Prop:monotony2} it follows
\beqq \max_{\vol(K)=c}\;\min_{q\in F^{cp}(K)}\ell_T(q) < \max_{\vol(K)=\vol(K^*)}\;\min_{q\in F^{cp}(K)}\ell_T(q)=\alpha.\eeqq

Conversely, let $K^*$ be a maximizer of
\beqq \max_{\vol(K)=c}\;\min_{q\in F^{cp}(K)}\ell_T(q),\eeqq
where $\alpha >0$ is chosen such that
\beqq \min_{q\in F^{cp}(K^*)}\ell_T(q)=:\widetilde{\alpha}<\alpha.\eeqq
Then we know from the above reasoning that $K^*$ is a minimizer of
\beqq \min_{T\in A(T,\widetilde{\alpha})}\vol(K),\eeqq
and from Proposition \ref{Prop:monotony1} it follows
\beqq \min_{K\in A(T,\alpha)}\vol(K)>\min_{K\in A(T,\widetilde{\alpha})}\vol(K)= \vol(K^*)=c.\eeqq
\epf

Hereinafter we will deal with the following two minimax problems\footnote{Whenever we write \beqq \min_{\vol(T)=d} \; \min_{K\in A(T,\alpha)} \vol(K)\eeqq or \beqq \max_{\vol(T)=d}\; \max_{\vol(K)=c} \; \min_{q\in F^{cp}(K)} \ell_T(q)\eeqq the minimum/maximum is understood to consider only convex bodies $T\subset\R^n$. This is implicitly indicated by the fact that we defined $A(\cdot,\alpha)$ and $\ell_{\cdot}(q)$ only for convex bodies $T\subset\R^n$.}: For $\alpha,d >0$ we will consider
\beqq \min_{\vol(T)=d} \; \min_{K\in A(T,\alpha)} \vol(K),\eeqq
and for $c,d >0$ we will consider
\beqq \max_{\vol(T)=d}\; \max_{\vol(K)=c} \; \min_{q\in F^{cp}(K)} \ell_T(q).\eeqq
It is indeed justified to write \say{$\min$} and \say{$\max$} respectively:

\bprop\label{Prop:minimaxproblem2}
Let $\alpha,d > 0$. Then we have
\beqq \inf_{\vol(T)=d} \; \min_{K\in A(T,\alpha)} \vol(K) = \min_{\vol(T)=d} \; \min_{K\in A(T,\alpha)} \vol(K).\eeqq
\eprop

\bpf
Let $(T_k)_{k\in\N}$ be a minimizing sequence of
\beq \inf_{\vol(T)=d} \; \min_{K\in A(T,\alpha)} \vol(K).\label{eq:minimaxproblem21}\eeq
Then there is an $R>0$ and a $k_0\in\N$ such that
\beqq T_k \subset B_R^n\quad \forall k\geq k_0.\eeqq

Indeed, if this is not the case, then there is a subsequence $\left(T_{k_j}\right)_{j\in\N}$ such that
\beq R_j:=\max\left\{R>0 : T_{k_j} \in F\left(B_R^n\right)\right\}\; \rightarrow \infty \quad (j\rightarrow \infty).\label{eq:minimaxproblem22}\eeq
This implies
\begin{align*}
V_j:&=\min\left\{\vol(K):K\in A\left(T_{k_j},\alpha\right)\right\}\\
&= \min\left\{\vol(K):K\in\mathcal{C}(\R^n),\; L_{T_{k_j}}(\alpha)\subseteq C(K)\right\}\\
&\rightarrow \infty \quad (j\rightarrow \infty).
\end{align*}
The latter follows from the fact that--\eqref{eq:minimaxproblem22} together with the convexity of $T_{k_j}$ for all $j\in\N$ and the constraint
\beqq \vol\left(T_{k_j}\right)=d\quad \forall j\in\N\eeqq
means that there are directions from the origin in which $T_{k_j}$ has to shrink for $j\rightarrow\infty$--for every $j\in\N$ we can find
\beqq q_j\in L_{T_{k_j}}(\alpha)\eeqq
($q_j$ can be constructed by using the aforementioned directions) for which
\beqq \ell_{T_{k_j}}(q_j)=\alpha\eeqq
means
\beqq \max_{t\in [0,\widetilde{T}_j]} |q_j(t)|\rightarrow\infty \quad (j\rightarrow\infty),\eeqq
where by $\widetilde{T}_j$ we denote the period of the closed curve $q_j$, and for every convex body $K_j\subset\R^n$ minimizing
\beqq \min\left\{\vol(K):K\in\mathcal{C}(\R^n),\; L_{T_{k_j}}(\alpha)\subseteq C(K)\right\}\eeqq
means
\beqq V_j=\vol(K_j)\rightarrow \infty \quad (j\rightarrow\infty).\eeqq
But this is not possible since $(T_k)_{k\in\N}$ is a minimizing sequence of \eqref{eq:minimaxproblem21}.

Then, we can apply Theorem \ref{Thm:Blaschkeselection}: There is a subsequence $\left(T_{k_l}\right)_{l\in\N}$ and a convex body $T\subset\R^n$ such that $T_{k_l}$ $d_H$-converges to $T$ for $l\rightarrow\infty$. We clearly have
\beqq \vol(T)=\vol\left(\lim_{l\rightarrow\infty}T_{k_l}\right)=\lim_{l\rightarrow\infty}\vol\left(T_{k_l}\right)=d.\eeqq
Therefore, $T$ is a minimizer of \eqref{eq:minimaxproblem21}.
\epf

\bprop\label{Prop:minimaxproblem3}
Let $c,d >0$. Then we have
\beq \sup_{\vol(T)=d}\; \max_{\vol(K)=c} \; \min_{q\in F^{cp}(K)} \ell_T(q) = \max_{\vol(T)=d}\; \max_{\vol(K)=c} \; \min_{q\in F^{cp}(K)} \ell_T(q).\label{eq:minimaxproblem30}\eeq
\eprop

\bpf
Let $\alpha >0$ and let us consider the minimax problem
\beq \min_{\vol(T)=d} \; \min_{K\in A(T,\alpha)} \vol(K).\label{eq:minimaxproblem31}\eeq
Let the pair
\beqq (K^*,T^*) \; \text{ with } \; \vol(T^*)=d \; \text{ and } \; K^*\in A(T^*,\alpha)\eeqq
be a minimizer of \eqref{eq:minimaxproblem31}, i.e., it is
\beqq \min_{\vol(T)=d}\; \min_{K\in A(T,\alpha)} \vol(K)=\min_{K\in A(T^*,\alpha)}\vol(K)=\vol(K^*)=:\widetilde{c}.\eeqq
By Theorem \ref{Thm:MainProperty1} $K^*$ is a maximizer of
\beqq \max_{\vol(K)=\widetilde{c}}\;\min_{q\in F^{cp}(K)}\ell_{T^*}(q)\eeqq
with
\beqq \min_{q\in F^{cp}(K^*)}\ell_{T^*}(q)=\alpha.\eeqq
Then, due to
\beqq \vol(T^*)=d\eeqq
we clearly have
\beq \alpha=\max_{\vol(K)=\widetilde{c}}\;\min_{q\in F^{cp}(K)}\ell_{T^*}(q)\leq \sup_{\vol(T)=d}\; \max_{\vol(K)=\widetilde{c}}\;\min_{q\in F^{cp}(K)}\ell_{T}(q).\label{eq:minimaxproblem32}\eeq
If this is a strict inequality, then there is a pair of convex bodies
\beqq (K^{**},T^{**}) \; \text{ with } \; \vol(T^{**})=d \; \text{ and } \; \vol(K^{**})=\widetilde{c}\eeqq
such that
\beqq \alpha < \max_{\vol(K)=\widetilde{c}}\;\min_{q\in F^{cp}(K)}\ell_{T^{**}}(q)=\min_{q\in F^{cp}(K^{**})}\ell_{T^{**}}(q)=:\widetilde{\alpha}.\eeqq
Then, by Theorem \ref{Thm:MainProperty1} $K^{**}$ is a minimizer of
\beqq \min_{K\in A(T^{**},\widetilde{\alpha})}\vol(K)\eeqq
with
\beqq \min_{K\in A(T^{**},\widetilde{\alpha})}\vol(K)=\vol(K^{**})=\widetilde{c}.\eeqq
Now, $\widetilde{\alpha} >\alpha$ together with Proposition \ref{Prop:monotony1} implies
\begin{align*}
\widetilde{c}=\vol(K^{**})=\min_{K\in A(T^{**},\widetilde{\alpha})}\vol(K)&\geq \min_{\vol(T)=d}\;\min_{K\in A(T,\widetilde{\alpha})}\vol(K)\\
&>\min_{\vol(T)=d}\;\min_{K\in A(T,\alpha)}\vol(K)\\
&=\min_{K\in A(T^*,\alpha)}\vol(K)\\
&=\vol(K^*)\\
&=\widetilde{c},
\end{align*}
a contradiction. Therefore, it follows that the inequality in \eqref{eq:minimaxproblem32} is in fact an equality, i.e.,
\beqq \sup_{\vol(T)=d}\; \max_{\vol(K)=\widetilde{c}}\;\min_{q\in F^{cp}(K)}\ell_{T}(q)=\alpha=\min_{q\in F^{cp}(K^*)}\ell_{T^*}(q).\eeqq
This means that the pair $(K^*,T^*)$ is a maximizer of
\beqq \sup_{\vol(T)=d}\; \max_{\vol(K)=\widetilde{c}} \; \min_{q\in F^{cp}(K)} \ell_T(q).\eeqq
Since it is sufficient to prove the claim \eqref{eq:minimaxproblem30} for one $c>0$, we are done.
\epf

\bthm\label{Thm:MainProperty2}
If the pair $(K^*,T^*)$ is a minimizer of
\beq \min_{\vol(T)=d}\; \min_{K\in A(T,\alpha)} \vol(K)\label{eq:minmin}\eeq
for $\alpha,d >0$, then $(K^*,T^*)$ is a maximizer of
\beq \max_{\vol(T)=d}\; \max_{\vol(K)=c}\; \min_{q\in F^{cp}(K)} \ell_T(q)\label{eq:maxmaxmin}\eeq
for
\beqq c:=\vol(K^*)\eeqq
with
\beqq \min_{q\in F^{cp}(K^*)}\ell_{T^*}(q)=\alpha.\eeqq

Conversely, if the pair $(K^*,T^*)$ is a maximizer of \eqref{eq:maxmaxmin} for $c,d>0$, then $(K^*,T^*)$ is a minimizer of \eqref{eq:minmin} for
\beqq \alpha:=\min_{q\in F^{cp}(K^*)}\ell_{T^*}(q)\eeqq
with
\beqq \vol(K^*)=c.\eeqq

Consequently, for $\alpha,c,d>0$ we have the equivalence
\beqq \min_{\vol(T)=d}\;\min_{K\in A(T,\alpha)}\vol(K)= c \; \Leftrightarrow \;  \max_{\vol(T)=d}\;\max_{\vol(K)=c}\;\min_{q\in F^{cp}(K)}\ell_T(q)= \alpha\eeqq
and moreover
\beq \min_{\vol(T)=d}\;\min_{K\in A(T,\alpha)}\vol(K)\geq c \; \Leftrightarrow \; \max_{\vol(T)=d}\;\max_{\vol(K)=c}\;\min_{q\in F^{cp}(K)}\ell_T(q)\leq \alpha.\label{eq:mainproperty21}\eeq
\ethm

\bpf
Let the pair $(K^*,T^*)$ be a minimizer of \eqref{eq:minmin} for $\alpha,d>0$, i.e., it is
\beqq \vol(T^*)=d \; \text{ and } \; K^* \in A(T^*,\alpha)\eeqq
such that
\beqq \min_{\vol(T)=d}\; \min_{K\in A(T,\alpha)} \vol(K)=\min_{K\in A(T^*,\alpha)}\vol(K)=\vol(K^*).\eeqq
Then, in the proof of Proposition \ref{Prop:minimaxproblem3} we have seen that $(K^*,T^*)$ is a maximizer of \eqref{eq:maxmaxmin} for
\beqq c:=\vol(K^*) \; \text{ with } \; \min_{q\in F^{cp}(K^*)}\ell_{T^*}(q)=\alpha.\eeqq


Conversely, let the pair $(K^*,T^*)$ be a maximizer of \eqref{eq:maxmaxmin} for $c,d>0$, i.e., $K^*,T^*\subset\R^n$ are convex bodies of volume $c$ and $d$, respectively, such that
\begin{align*}
\max_{\vol(T)=d}\; \max_{\vol(K)=c}\; \min_{q\in F^{cp}(K)} \ell_T(q)&=\max_{\vol(K)=c}\;\min_{q\in F^{cp}(K)}\ell_{T^*}(q)\\
&=\min_{q\in F^{cp}(K^*)}\ell_{T^*}(q)\\
&=:\alpha.
\end{align*}
By Theorem \ref{Thm:MainProperty1} $K^*$ minimizes
\beqq \min_{K\in A(T^*,\alpha)}\vol(K)\eeqq
with
\beqq \vol(K^*)=c.\eeqq
Then, we clearly have
\beqq c=\vol(K^*)=\min_{K\in A(T^*,\alpha)}\vol(K)\geq \min_{\vol(T)=d}\;\min_{K\in A(T,\alpha)}\vol(K).\eeqq
If this is a strict inequality, then there is a pair $(K^{**},T^{**})$ with
\beqq c>\min_{\vol(T)=d}\;\min_{K\in A(T,\alpha)}\vol(K)=\min_{K\in A(T^{**},\alpha)}\vol(K)=\vol(K^{**})=:\widetilde{c},\eeqq
where
\beqq K^{**}\in A(T,\alpha)\eeqq
and $T^{**}\subset\R^n$ is a convex body of volume $d$. Then, by Theorem \ref{Thm:MainProperty1} $K^{**}$ is a maximizer of
\beqq \max_{\vol(K)=\widetilde{c}}\;\min_{q\in F^{cp}(K)}\ell_{T^{**}}(q)\eeqq
with
\beqq \max_{\vol(K)=\widetilde{c}}\;\min_{q\in F^{cp}(K)}\ell_{T^{**}}(q)=\min_{q\in F^{cp}(K^{**})}\ell_{T^{**}}(q)=\alpha.\eeqq
Now, $\widetilde{c}<c$ together with Proposition \ref{Prop:monotony2} implies
{\allowdisplaybreaks\begin{align*}
\alpha=\min_{q\in F^{cp}(K^{**})}\ell_{T^{**}}(q)&=\max_{\vol(K)=\widetilde{c}}\;\min_{q\in F^{cp}(K)}\ell_{T^{**}}(q)\\
&\leq \max_{\vol(T)=d}\;\max_{\vol(K)=\widetilde{c}}\;\min_{q\in F^{cp}(K)}\ell_{T}(q)\\
&<\max_{\vol(T)=d}\; \max_{\vol(K)=c}\; \min_{q\in F^{cp}(K)} \ell_T(q)\\
&=\max_{\vol(K)=c}\;\min_{q\in F^{cp}(K)}\ell_{T^*}(q)\\
&=\min_{q\in F^{cp}(K^*)}\ell_{T^*}(q)\\
&=\alpha,
\end{align*}}%
a contradiction. Therefore,
\beqq \min_{\vol(T)=d}\;\min_{K\in A(T,\alpha)}\vol(K)=c=\vol(K^*)=\min_{K\in A(T^*,\alpha)}\vol(K),\eeqq
i.e., the pair $(K^*,T^*)$ is a minimizer of \eqref{eq:minmin}.

From the before proven it clearly follows the equivalence
\beqq \min_{\vol(T)=d}\;\min_{K\in A(T,\alpha)}\vol(K)= c \;\Leftrightarrow\; \max_{\vol(T)=d}\;\max_{\vol(K)=c}\;\min_{q\in F^{cp}(K)}\ell_T(q)= \alpha.\eeqq
for $\alpha, c, d >0$.

In order to prove \eqref{eq:mainproperty21} it is sufficient to show
\beqq \min_{\vol(T)=d}\;\min_{K\in A(T,\alpha)}\vol(K)> c \;\Leftrightarrow\; \max_{\vol(T)=d}\;\max_{\vol(K)=c}\;\min_{q\in F^{cp}(K)}\ell_T(q)< \alpha.\eeqq

Let the pair $(K^*,T^*)$ be a minimizer of
\beqq \min_{\vol(T)=d}\;\min_{K\in A(T,\alpha)}\vol(K),\eeqq
where $c>0$ is chosen such that
\beqq c<\min_{\vol(T)=d}\;\min_{K\in A(T,\alpha)}\vol(K)=\min_{K\in A(T^*,\alpha)}\vol(K)=\vol(K^*)=:\widetilde{c}.\eeqq
From above reasoning we know that $(K^*,T^*)$ maximizes \eqref{eq:maxmaxmin} (for $c$ replaced by $\widetilde{c}$), i.e., $K^*,T^*\subset\R^n$ are convex bodies of volume $\widetilde{c}$ and $d$, respectively, such that
\begin{align*}
\max_{\vol(T)=d}\;\max_{\vol(K)=\widetilde{c}}\;\min_{q\in F^{cp}(K)}\ell_T(q)&=\max_{\vol(K)=\widetilde{c}}\;\min_{q\in F^{cp}(K)}\ell_{T^*}(q)\\
&=\min_{q\in F^{cp}(K^*)}\ell_{T^*}(q)\\
&=\alpha.
\end{align*}
Now, $c<\widetilde{c}$ together with Proposition \ref{Prop:monotony2} implies
\beqq \max_{\vol(T)=d}\;\max_{\vol(K)=c}\;\min_{q\in F^{cp}(K)}\ell_T(q) < \max_{\vol(T)=d}\;\max_{\vol(K)=\widetilde{c}}\;\min_{q\in F^{cp}(K)}\ell_T(q)=\alpha.\eeqq

Conversely, let $(K^*,T^*)$ be a maximizer of
\beqq \max_{\vol(T)=d}\;\max_{\vol(K)=c}\;\min_{q\in F^{cp}(K)}\ell_T(q),\eeqq
i.e., $K^*,T^*\subset\R^n$ are convex bodies of volume $c$ and $d$, respectively, where $\alpha >0$ is chosen such that
\begin{align*}
\alpha > \max_{\vol(T)=d}\;\max_{\vol(K)=c}\;\min_{q\in F^{cp}(K)}\ell_T(q)&=\max_{\vol(K)=c}\;\min_{q\in F^{cp}(K)}\ell_{T^*}(q)\\
&=\min_{q\in F^{cp}(K^*)}\ell_{T^*}(q)\\
&=:\widetilde{\alpha}.
\end{align*}
Then we know from above reasoning that $(K^*,T^*)$ minimizes \eqref{eq:minmin} (for $\alpha$ replaced by $\widetilde{\alpha}$), i.e.,
\beqq \min_{\vol(T)=d}\;\min_{K\in A(T,\widetilde{\alpha})}\vol(K)=\min_{K\in A(T^*,\widetilde{\alpha})}\vol(K)=\vol(K^*)=c.\eeqq
Now, $\alpha > \widetilde{\alpha}$ together with Proposition \ref{Prop:monotony1} implies
\beqq \min_{\vol(T)=d}\;\min_{K\in A(T,\alpha)}\vol(K)>\min_{\vol(T)=d}\;\min_{K\in A(T,\widetilde{\alpha})}\vol(T)=c.\eeqq
\epf

\section{Proofs of Theorems \ref{Thm:main1}, \ref{Thm:Mahler}, \ref{Thm:relations} and Corollary \ref{Cor:optimizationmain1}}\label{Sec:main1}

In the following, we mainly make use of Theorems \ref{Thm:MainProperty1} and \ref{Thm:MainProperty2}. However, we begin by rewriting Viterbo's conjecture for convex Lagrangian products:

\bprop\label{Prop:rewriteViterbo}
Viterbo's conjecture for convex Lagrangian products $K\times T \subset \R^n\times\R^n$
\beqq \vol(K\times T) \geq \frac{c_{EHZ}(K\times T)^n}{n!},\quad K,T \in \mathcal{C}(\R^n),\eeqq
is equivalent to
\beqq \max_{\vol(K)=1}\; \max_{\vol(T)=1}\; \min_{q\in F^{cp}(K)} \ell_T(q) \leq \sqrt[n]{n!},\quad K,T \in \mathcal{C}(\R^n).\eeqq
\eprop

\bpf
Using Proposition \ref{Prop:invViterbo}, Viterbo's conjecture for convex Lagrangian products is equivalent to
\beqq \max_{\vol(K)=1}\; \max_{\vol(T)=1}\; c_{EHZ}(K\times T) \leq \sqrt[n]{n!},\quad K,T \in \mathcal{C}(\R^n).\eeqq
By Theorem \ref{Thm:relationcapacity}, this is equivalent to
\beqq \max_{\vol(K)=1}\; \max_{\vol(T)=1}\; \min_{q\in F^{cp}(K)} \ell_T(q) \leq \sqrt[n]{n!},\quad K,T \in \mathcal{C}(\R^n).\eeqq
\epf

\bpf[Proof of Theorem \ref{Thm:main1}]
Using Proposition \ref{Prop:rewriteViterbo}, Viterbo's conjecture for convex Lagrangian products is equivalent to
\beqq \max_{\vol(K)=1}\; \max_{\vol(T)=1}\; \min_{q\in F^{cp}(K)} \ell_T(q) \leq \sqrt[n]{n!},\quad K,T \in \mathcal{C}(\R^n).\eeqq
After applying Theorem \ref{Thm:MainProperty2}, it is further equivalent to
\beqq \min_{\vol(T)=1}\;\min_{K\in A\left(T,\sqrt[n]{n!}\right)}\vol(K)\geq 1,\quad K,T \in \mathcal{C}(\R^n).\eeqq
Using Proposition \ref{Prop:Ahomogenity2}, this can be written as
\beqq \min_{\vol(T)=1}\;\min_{K\in A(T,1)}\vol(K)\geq \frac{1}{n!},\quad K,T \in \mathcal{C}(\R^n).\eeqq

By similar reasoning, Theorem \ref{Thm:MainProperty2} also guarantees the equivalence of the equality case of Viterbo's conjecture for convex Lagrangian products $K\times T\subset\R^n\times\R^n$
\beq \vol(K\times T) = \frac{c_{EHZ}(K\times T)^n}{n!},\quad K,T \in \mathcal{C}(\R^n),\label{eq:rewriteViterbo1}\eeq
i.e.,
\beqq \max_{\vol(K)=1}\; \max_{\vol(T)=1}\; \min_{q\in F^{cp}(K)} \ell_T(q) = \sqrt[n]{n!},\quad K,T \in \mathcal{C}(\R^n),\eeqq
and
\beq \min_{\vol(T)=1}\;\min_{K\in A(T,1)}\vol(K)= \frac{1}{n!},\quad K,T \in \mathcal{C}(\R^n).\label{eq:rewriteViterbo2}\eeq

Moreover, Theorem \ref{Thm:MainProperty2} guarantees the following: If $K^*\times T^*$ is a solution of \eqref{eq:rewriteViterbo1} satisfying
\beq \vol(K^*)=\vol(T^*)=1\label{eq:rewriteViterbo3}\eeq
(note that, applying Proposition \ref{Prop:invViterbo}, the property of being a solution of \eqref{eq:rewriteViterbo1} is invariant under scaling), then the pair $(K^*,T^*)$ is a solution of \eqref{eq:rewriteViterbo2}. And conversely, if the pair $(K^*,T^*)$ is a solution of \eqref{eq:rewriteViterbo2}, then $K^*\times T^*$ is a solution of \eqref{eq:rewriteViterbo1}.
\epf

\bpf[Proof of Corollary \ref{Cor:optimizationmain1}]
In view of the proof of Theorem \ref{Thm:main1}, for convex bodies $K,T\subset\R^n$, it is sufficient to prove the following equality:
\beq \min_{\vol(T)=1}\; \min_{K\in A(T,1)} \vol(K) = \min_{\vol(T)=1}\;\min_{a_q\in\R^n}\vol\bigg(\conv\bigg\{ \bigcup_{q\in L_T(1)} (q+a_q)\bigg\}\bigg).\label{eq:approach0}\eeq
But this follows from the following gradually observation: First, we notice that the volume-minimizing convex cover for a set of closed curves is, equivalently, the volume-minimizing convex hull of these closed curves. So, if we ask for lower bounds of
\beqq \min_{K\in A(T,1)}\vol(K),\eeqq
we note that for $q_1,...,q_k\in L_T(1)$, we have
\beqq \min_{(a_1,...,a_k)\in (\R^n)^k}\;\vol\left(\conv\{q_1+a_1,...,q_k+a_k\}\right) \leq \min_{K\in A(T,1)}\vol(K).\eeqq
This estimate can be further improved by
\beqq \max_{q_1,...,q_k\in L_T(1)}\;\min_{(a_1,...,a_k)\in (\R^n)^k}\;\vol(\conv\{q_1+a_1,...,q_k+a_k\}) \leq \min_{K\in A(T,1)}\vol(K),\eeqq
so that eventually we get
\beqq \min_{a_q\in\R^n}\vol\bigg(\conv\bigg\{ \bigcup_{q\in L_T(1)} (q+a_q)\bigg\}\bigg)=\min_{K\in A(T,1)}\vol(K),\eeqq
where the minimum on the left runs for every $q\in L_T(1)$ over all possible translations in $\R^n$. Minimizing this equation over all convex bodies $T\subset\R^n$ of volume $1$, we get \eqref{eq:approach0}.
\epf

\bpf[Proof of Theorem \ref{Thm:Mahler}]
Because of
\beqq c_{EHZ}(T\times T^\circ)=4\eeqq
for all centrally symmetric convex bodies $T\subset\R^n$ (cf.\;\cite{ArtKarOst2013}), Mahler's conjecture for centrally symmetric convex bodies is equivalent to
\beq \vol(T\times T^\circ)\geq \frac{c_{EHZ}(T\times T^\circ)^n}{n!}, \quad T\in \mathcal{C}^{cs}(\R^n).\label{eq:mahlerpr1}\eeq
Fixing
\beqq \vol(T)=1,\eeqq
which is without loss of generality due to Proposition \ref{Prop:invMahler}, and using Theorem \ref{Thm:relationcapacity}, \eqref{eq:mahlerpr1} is equivalent to
\beqq \sqrt[n]{n! \vol(T^\circ)}\geq c_{EHZ}(T\times T^\circ)=\min_{q\in F^{cp}(T)}\ell_{T^\circ}(q),\quad T\in \mathcal{C}^{cs}(\R^n).\eeqq
This can be written as
\beqq \max_{\vol(T)=1}\;\min_{q\in F^{cp}(T)}\ell_{T^\circ}(q)\leq \sqrt[n]{n! \vol(T^\circ)},\quad T\in \mathcal{C}^{cs}(\R^n),\eeqq
which, by Theorem \ref{Thm:MainProperty1}, is equivalent to
\beqq \min_{T\in A\left(T^\circ,\sqrt[n]{n! \vol(T^\circ)}\right)}\vol(T)\geq 1,\quad T\in \mathcal{C}^{cs}(\R^n).\eeqq
Applying Proposition \ref{Prop:Ahomogenity2}, we finally conclude that Mahler's conjecture for centrally symmetric convex bodies is equivalent to
\beqq \min_{T\in A\left(T^\circ,\sqrt[n]{\vol(T^\circ)}\right)}\vol(T)\geq \frac{1}{n!},\quad T\in \mathcal{C}^{cs}(\R^n).\eeqq

By similar reasoning, Theorem \ref{Thm:MainProperty1} also guarantees the equivalence of the equality case of Mahler's conjecture for centrally symmetric convex bodies $T\subset\R^n$
\beq \vol(T)\vol(T^\circ) =\frac{4^n}{n!},\label{eq:mahlerpr2}\eeq
i.e.,
\beqq \max_{\vol(T)=1}\;\min_{q\in F^{cp}(T)}\ell_{T^\circ}(q)= \sqrt[n]{n! \vol(T^\circ)},\quad T\in \mathcal{C}^{cs}(\R^n),\eeqq
and
\beq \min_{T\in A(T^\circ,\sqrt[n]{\vol(T^\circ)})}\vol(T) = \frac{1}{n!},\quad T\in \mathcal{C}^{cs}(\R^n).\label{eq:mahlerpr3}\eeq

Moreover, Theorem \ref{Thm:MainProperty1} guarantees the following: If $T^*$ is a solution of \eqref{eq:mahlerpr2} satisfying
\beqq \vol(T^*)=1\eeqq
(note that, applying Proposition \ref{Prop:invMahler}, the property of being a solution of \eqref{eq:mahlerpr2} is invariant under scaling), then it is a solution of \eqref{eq:mahlerpr3}. And conversely, if $T^*$ is a solution of \eqref{eq:mahlerpr3}, then it is also a solution of \eqref{eq:mahlerpr2}.
\epf

\bpf[Proof of Theorem \ref{Thm:relations}]
The equivalence of (i), (ii), and (iii) follows from Theorem \ref{Thm:relationcapacity}. The equivalence of (i) and (iv) follows from Theorem \ref{Thm:MainProperty1}. The equivalence of (iv) and (v) can be concluded as within the proof of Corollary \ref{Cor:optimizationmain1}. For the case of strictly convex $T\subset\R^n$, the equivalence of (i) and (vi) follows from \cite[Theorem 1.3]{KruppRudolf2022}.

The addition that every equality case $(K^*,T^*)$ of any of the inequalities is also an equality case of all the others is guaranteed by Theorem \ref{Thm:MainProperty1}.
\epf

\section{Proof of Theorem \ref{Thm:main2}}\label{Sec:main2}

We recall a sligthly rephrased version of the main result of Haim-Kislev in \cite{Haim-Kislev2019}:

\bthm\label{Thm:HaimKislev}
Let $P\subset\R^{2n}$ be a convex polytope. Then, there is an action-minimizing closed characteristic $x$ on $\partial P$ which is a closed polygonal curve consisting of finitely many segments
\beqq [x(t_j),x(t_{j+1})]\eeqq
given by
\beqq x(t_{j+1})=x(t_j)+\lambda_j J\nabla H_P(x_j),\quad \lambda_j >0,\eeqq
while $x_j\in \mathring{F}_j$, $F_j$ is a facet of $P$ and $x$ visits every facet $F_j$ at most once.
\ethm

For the proof of Theorem \ref{Thm:main2}, we need the following theorem:

\bthm\label{Thm:Estimate}
If $P\subset\R^{2n}$ is a convex polytope, then we have
\beqq c_{EHZ}(P)= R_P \min_{q\in F^{cp}(P)}\ell_{\frac{JP}{2}}(q)\eeqq
with
\beqq R_P = \frac{\min_{q\in F^{cp}_*(P)}\ell_{\frac{JP}{2}}(q)}{\min_{q\in F^{cp}(P)}\ell_{\frac{JP}{2}}(q)}\geq 1.\eeqq
\ethm

We remark that, in light of Theorem \ref{Thm:relationcapacity}, Theorem \ref{Thm:Estimate} implies the following relationship between the EHZ-capacity of $P$ and the EHZ-capacity of the Lagrangian product $P\times \frac{1}{2}JP$:
\beqq c_{EHZ}(P)=R_P c_{EHZ}\left(P\times \frac{1}{2}JP\right).\eeqq

For the proof of Theorem \ref{Thm:Estimate}, we need the following proposition. We remark that in the proof of Theorem \ref{Thm:Estimate}, we need it only in the case of action-minimizing closed characteristics on the boundary of a polytope. However, we will state it in full generality which has relevance beyond its use in the proof of Theorem \ref{Thm:Estimate} (which we will briefly address below).

\bprop\label{Prop:characteristicgeodesic}
Let $C\subset\R^{2n}$ be a convex body. Let $x$ be any closed characteristic on $\partial C$. Then, the action of $x$ equals its $\ell_{\frac{JC}{2}}$-length:
\beqq \A(x)=\ell_{\frac{JC}{2}}(x).\eeqq
\eprop

Proposition \ref{Prop:characteristicgeodesic} implies a noteworthy connection between closed characteristics and closed Finsler geodesics: Every closed characteristic on $\partial C$ can be interpreted as a Finsler geodesic with respect to the Finsler metric determined by $\mu_{2JC^\circ}$ and which is parametrized by arc length. This raises a number of questions; for example, which Finsler geodesics are closed characteristics (we note that, usually, there are more geodesics than those which, by the least action principle and Proposition \ref{Prop:characteristicgeodesic}, can be associated to closed characteristics) and whether the length-minimizing Finsler geodesics are of this kind. Following this line of thought, would lead to the question whether it is possible to deduce Viterbo's conjecture from systolic inequalities for Finsler geodesics. However, we leave these questions for further research.

\bpf[Proof of Proposition \ref{Prop:characteristicgeodesic}]
By
\beqq \dot{x}(t)\in J\partial H_C(x(t))\quad \text{a.e.},\eeqq
we conclude
\beqq \frac{1}{2}\left(\mu_{2JC^\circ}(\dot{x}(t))\right)^2=H_{2JC^\circ}(\dot{x}(t))\in H_{2JC^\circ}(J\partial H_C(x(t)))=\frac{1}{4} H_{C^\circ}(\partial H_C(x(t)))\quad \text{a.e.},\eeqq
where we used the facts
\beqq J^{-1}=-J,\quad H_C(Jx)=H_{J^{-1}C}(x) \eeqq
and
\beqq H_{\lambda C}(x)=H_{C}\left(\frac{1}{\lambda}x\right)=\frac{1}{\lambda^2}H_C(x),\quad \lambda \neq 0,\eeqq
(cf.\;\cite[Proposition 2.3(iii)]{KruppRudolf2022}). From Proposition \ref{Prop:ConsEuler}, we therefore conclude
\beqq \frac{1}{2}\left(\mu_{2JC^\circ}(\dot{x}(t))\right)^2=\frac{1}{4}H_C(x(t))=\frac{1}{8} \quad \text{a.e.}\eeqq
and consequently
\beqq \mu_{2JC^\circ}(\dot{x}(t))=\frac{1}{2} \quad \text{a.e.}\eeqq
Considering
\beqq (2JC^\circ)^\circ = \frac{1}{2}JC\eeqq
(cf.\;\cite[Proposition 2.1]{KruppRudolf2022}), we obtain
\begin{align*}
\ell_{\frac{JC}{2}}(x)=\int_0^T \mu_{\left(\frac{JC}{2}\right)^\circ}(\dot{x}(t))\,\mathrm{d}t = \int_0^T \mu_{2JC^\circ}(\dot{x}(t))\,\mathrm{d}t=\int_0^T \frac{1}{2}\,\mathrm{d}t=\frac{T}{2}=\A(x),
\end{align*}
where the last equality follows from
\beqq \A(x)=-\frac{1}{2}\int_0^T \langle J\dot{x}(t),x(t)\rangle \, \mathrm{d}t \in \frac{1}{2}\int_0^T \langle \partial H_C(x(t)),x(t)\rangle\, \mathrm{d}t\eeqq
which by Proposition \ref{Prop:Euleridentity} and the $2$-homogeneity of $H_C$ implies
\beqq  \A(x) =\int_0^T H_C(x(t)) \,\mathrm{d}t = \frac{T}{2}.\eeqq
\epf

Then, we come to the proof of Theorem \ref{Thm:Estimate}:

\bpf[Proof of Theorem \ref{Thm:Estimate}]
The idea behind the proof is to associate action-minimizing closed characteristics on $\partial P$ in the sense of Theorem \ref{Thm:HaimKislev} with $\ell_{\frac{1}{2}JP}$-minimizing closed $(P,\frac{JP}{2})$-Minkowski billiard trajectories.

Let $x$ be an action-minimizing closed characteristic on $\partial P$ in the sense of Theorem \ref{Thm:HaimKislev}. Let us assume $x$ is moving on the facets of $P$ according to the order
\beqq F_1\rightarrow F_2 \rightarrow ... \rightarrow F_m \rightarrow F_1,\eeqq
while the linear flow on every facet is given by the $J$-rotated normal vector at the interior of this facet. Out of every trajectory segment
\beqq \orb(x)\cap \mathring{F}_j,\eeqq
we choose one point $q_j$ arbitrarily (on the whole requiring $q_i\neq q_j$ for $i\neq j$) and connect these points by straight lines (by maintaining the order of the corresponding facets) constructing a closed polygonal curve
\beqq q:=(q_1,...,q_m)\eeqq
within $P$ which has its vertices on $\partial P$. From Lemma \ref{Lem:Pre2} (which we provide subsequently), we derive
\beqq \ell_{\frac{JP}{2}}(q)=\ell_{\frac{JP}{2}}(x)\eeqq
since the trajectory segment of $x$ between the two consecutive points $q_j$ and $q_{j+1}$--let us call it $\orb(x)_{q_j\rightarrow q_{j+1}}$--together with the line from $q_j$ to $q_{j+1}$ (as trajectory segment of $q$)--let us call it $[q_j,q_{j+1}]$--builds a triangle with the property that
\beqq \mu_{2JP^\circ}\left(\orb(x)_{q_j\rightarrow q_{j+1}}\right)=\mu_{2JP^\circ}([q_j,q_{j+1}]).\eeqq   
We therefore conclude from Proposition \ref{Prop:characteristicgeodesic} that
\beqq \ell_{\frac{JP}{2}}(q) = \A(x).\eeqq
Because of the arbitrariness of the choice of $q_j$ within $\orb(x)\cap \mathring{F}_j$, we can assign infinitely many different closed polygonal curves of the above kind to one action-minimizing closed characteristic fulfilling the demanded conditions.

Each of these closed polygonal curves $q$ is a closed $(P\times \frac{1}{2}JP)$-Minkowski billiard trajectory: This follows from the fact that $q$ fulfills
\beqq \begin{cases}q_{j+1}-q_j\in N_{\frac{1}{2}JP}(p_j),\\p_{j+1}-p_j\in -N_P(q_{j+1}),\end{cases}\eeqq
for the closed polygonal curve $p=(p_1,...,p_m)$ in $\frac{1}{2}JP$ with
\beqq p_{j-1}\in \partial \left(\frac{1}{2}JP\right)\eeqq
given as the intersection point
\beqq \frac{1}{2}J\left(\{q_{j-1}+tJ\nabla H_P(q_{j-1}):t\in\R\}\cap \{q_{j}+tJ\nabla H_P(q_{j}):t\in\R\}\right)  \subset \frac{JF_{j-1}}{2} \cap \frac{JF_{j}}{2} \eeqq
for all $j\in\{2,...,m+1\}$.

Indeed, from the definition of $p$, it follows
\beq p_{j+1}-p_j\in -N_P(q_{j+1}) \quad \forall j\in\{1,...,m\}\label{eq:minusexplanation}\eeq
since, by construction, $p_{j+1}-p_j$ is a multiple of the outer normal vector at $P$ in $q_j$ rotated by twofold multiplication with $J$ ($J^2=-\mathbb{1}$ produces the minus sign in \eqref{eq:minusexplanation}). Since, by construction,
\beqq J^{-1}(q_j-q_{j-1})\eeqq
is in the normal cone at $P$ in the intersection point
\beqq \{q_{j-1}+tJ\nabla H_P(q_{j-1}):t\in\R\}\cap \{q_{j}+tJ\nabla H_P(q_{j}):t\in\R\} \subset F_{j-1} \cap F_{j},\eeqq
roation by $\frac{1}{2}J$ then implies that $q_j-q_{j-1}$ is in the normal cone at $\frac{1}{2}JP$ in $p_{j-1}$. This implies 
\beqq q_{j}-q_{j-1}\in N_P(p_{j-1})\quad \forall j\in\{1,...,m\}.\eeqq

\begin{figure}[h!]
\centering
\def\svgwidth{420pt}
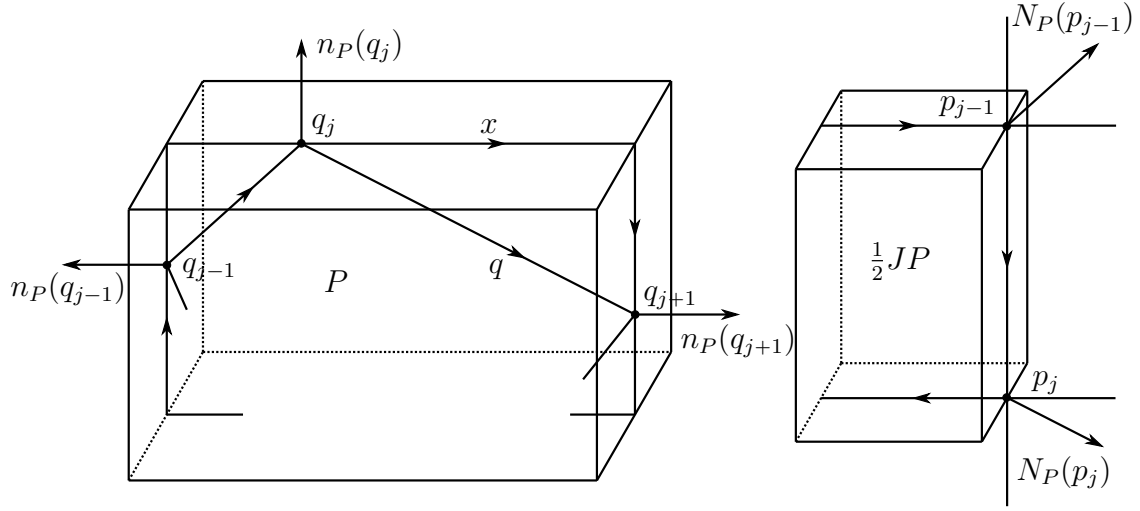
\caption[Association of closed characteristics to closed Minkowski billiard trajectories]{$q=(q_1,...,q_m)$ is a closed $(P,\frac{1}{2}JP)$-Minkowski billiard trajectory with $p=(p_1,...,p_m)$ as its dual billiard trajectory in $\frac{1}{2}JP$.}
\label{img:Billiard}
\end{figure}

From \cite[Proposition 3.9]{KruppRudolf2022}, it follows that $q$ cannot be translated into $\mathring{P}$, i.e.,
\beqq q\in F^{cp}(P).\eeqq
From the construction of $q$, we moreover know
\beq q\in F^{cp}_*(P),\label{eq:wormoneinequality}\eeq
where we recall that $F^{cp}_*(P)$ as subset of $F^{cp}(P)$ was defined as the set of all closed polygonal curves $q=(q_1,...,q_m)$ in $F^{cp}(P)$ for which $q_j$ and $q_{j+1}$ are on neighbouring facets $F_j$ and $F_{j+1}$ of $P$ such that there are $\lambda_j,\mu_{j+1}\geq 0$ with
\beqq q_{j+1}=q_j + \lambda_j J\nabla H_P(x_j)+\mu_{j+1}J\nabla H_P(x_{j+1}),\eeqq
where $x_j$ and $x_{j+1}$ are arbitrarily chosen interior points of $F_j$ and $F_{j+1}$, respectively.

Because of \eqref{eq:wormoneinequality}, we have 
\beqq \ell_{\frac{JP}{2}}(q) \geq \min_{\widetilde{q}\in F^{cp}_*(P)}\ell_{\frac{JP}{2}}(\widetilde{q}).\eeqq
Since, by definition and the above considerations, every closed polygonal curve in $F^{cp}_*(P)$ is associated with a closed characteristic on $\partial P$, where the $\ell_{\frac{JP}{2}}$-length of the former coincides with the action of the latter, and $x$ (to which $q$ is associated) was chosen to be action-minimizing, we actually have
\beqq \ell_{\frac{JP}{2}}(q) = \min_{\widetilde{q}\in F^{cp}_*(P)}\ell_{\frac{JP}{2}}(\widetilde{q}).\eeqq

Altogether, this implies
\beqq c_{EHZ}(P)=\A(x)=\ell_{\frac{JP}{2}}(x)=\ell_{\frac{JP}{2}}(q)=\min_{\widetilde{q}\in F^{cp}_*(P)}\ell_{\frac{JP}{2}}(\widetilde{q})=R_P \min_{\widetilde{q}\in F^{cp}(P)}\ell_{\frac{JP}{2}}(\widetilde{q})\eeqq
for
\beqq R_P = \frac{\min_{q\in F^{cp}_*(P)}\ell_{\frac{JP}{2}}(q)}{\min_{q\in F^{cp}(P)}\ell_{\frac{JP}{2}}(q)}\geq 1.\eeqq
\epf

\blem\label{Lem:Pre2}
Let $P\subset\R^{2n}$ be a convex polytope. If
\beqq y=\lambda_iJ\nabla H_P(x_i)+\lambda_j J\nabla H_P(x_j), \; \lambda_i,\lambda_j\geq 0,\eeqq
where $F_i$ and $F_j$ are neighbouring facets of $P$ with $x_i\in \mathring{F}_i$ and $x_j\in \mathring{F}_j$, then
\beqq \mu_{2JP^{\circ}}(y)=\lambda_i\mu_{2JP^{\circ}}(J\nabla H_P(x_i))+\lambda_j\mu_{2JP^{\circ}}(J\nabla H_P(x_j))=\frac{1}{2}(\lambda_i + \lambda_j).\eeqq
\elem

\bpf
We first notice that
\beqq \nabla H_P(x_i) \; \text{ and } \; \nabla H_P(x_j)\eeqq
are neighbouring vertices of $P^\circ$, i.e.,
\beqq t \nabla H_P(x_i)+(1-t)\nabla H_P(x_j)\in \partial P^\circ\;\;\; \forall t\in[0,1].\eeqq

Indeed, from the fact that $\nabla H_P(x_i)$ and $\nabla H_P(x_j)$ are elements of the one dimensional normal cone at $\mathring{F}_i$ and $\mathring{F}_j$, we conclude by the properties of the polar of convex polytopes (cf.\;\cite[Chapter 3.3]{GallierQuaintance2017}) that they point into the direction of two neigbouring vertices of $P^\circ$. Using Proposition \ref{Prop:ConsEuler}, we calculate
\beqq H_{P^\circ}(\nabla H_P(x_{i}))=H_P(x_{i})=\frac{1}{2}\eeqq
and
\beqq H_{P^\circ}(\nabla H_P(x_{j}))=H_P(x_{j})=\frac{1}{2}\eeqq
and conclude that $\nabla H_P(x_i)$ and $\nabla H_P(x_j)$ actually are these two neighbouring vertices of $P^\circ$. 

Using for convex body $C\subset\R^{2n}$ and $\lambda > 0$ the properties
\beqq \mu_{\lambda C}(x)=\frac{1}{\lambda}\mu_{C}(x) \;\text{ and }\;\mu_{J C}(Jx)=\mu_C(x),\;\; x\in\R^{2n},\eeqq
(cf.\;\cite[Proposition 2.3(iii)]{KruppRudolf2022}), we derive
\allowdisplaybreaks{\begin{align*}
\mu_{2JP^\circ}(y)&=\mu_{2JP^\circ}(\lambda_iJ\nabla H_P(x_i)+\lambda_j J\nabla H_P(x_j))\\
&= \mu_{2P^\circ}(\lambda_i\nabla H_P(x_i)+\lambda_j\nabla H_P(x_j))\\
&=\frac{1}{2}\left(\mu_{P^\circ}(\lambda_i\nabla H_P(x_i)+\lambda_j\nabla H_P(x_j))\right)\\
&\stackrel{(\star)}{=} \frac{1}{2}\left(\mu_{P^{\circ}}(\lambda_i\nabla H_P(x_i))+\mu_{P^{\circ}}(\lambda_j\nabla H_P(x_j))\right)\\
&= \frac{1}{2}\left(\lambda_i\mu_{P^{\circ}}(\nabla H_P(x_i))+\lambda_j\mu_{P^{\circ}}(\nabla H_P(x_j))\right)\\
&=\frac{1}{2}(\lambda_i+\lambda_j),
\end{align*}}%
where in $(\star)$ we used that, by the choice of $x_i$ and $x_j$ and the properties of polar bodies, $\nabla H_P(x_i)$ and $\nabla H_P(x_j)$ are neighbouring vertices of $P^\circ$ and, therefore, in $(\star)$, the initial term can be splitted linearly.
\epf

\bpf[Proof of Theorem \ref{Thm:main2}]
Viterbo's conjecture for convex polytopes in $\R^{2n}$ can be written as
\beqq \vol(P)\geq \frac{c_{EHZ}(P)^n}{n!},\quad P\in\mathcal{C}^p\left(\R^{2n}\right),\eeqq
which, by Theorem \ref{Thm:Estimate}, is equivalent to
\beqq \vol(P)\geq \frac{R_P^n}{2^n n!}c_{EHZ}(P\times JP)^n,\quad P\in\mathcal{C}^p\left(\R^{2n}\right).\eeqq
By referring to Proposition \ref{Prop:invViterbo}, we can assume
\beqq \vol(P)=1\eeqq
without loss of generality and get
\beqq c_{EHZ}(P\times JP)\leq \frac{2\sqrt[n]{n!}}{R_P},\quad P\in\mathcal{C}^p\left(\R^{2n}\right),\eeqq
which, by Theorem \ref{Thm:relationcapacity}, is equivalent to
\beqq \max_{\vol(P)=1}\; \min_{q\in F^{cp}(P)}\ell_{JP}(q) \leq \frac{2\sqrt[n]{n!}}{R_P},\quad P\in\mathcal{C}^p\left(\R^{2n}\right).\eeqq
By Theorem \ref{Thm:MainProperty1}, this is equivalent to
\beqq \min_{P\in A\left(JP,\frac{2\sqrt[n]{n!}}{R_P}\right)}\vol(P)\geq 1,\quad P\in\mathcal{C}^p\left(\R^{2n}\right),\eeqq
and, after applying Proposition \ref{Prop:Ahomogenity2}, to
\beqq \min_{P\in A\left(JP,\frac{1}{R_{P}}\right)}\vol(P)\geq \frac{1}{2^n n!},\quad P\in\mathcal{C}^p\left(\R^{2n}\right).\eeqq

By similar reasoning, Theorem \ref{Thm:MainProperty1} also guarantees the equivalence of
\beq \max_{\vol(P)=1}\; \min_{q\in F^{cp}(P)}\ell_{JP}(q) = \frac{2\sqrt[n]{n!}}{R_P},\quad P\in\mathcal{C}^p\left(\R^{2n}\right),\label{eq:proofmain21}\eeq
and
\beq \min_{P\in A\left(JP,\frac{1}{R_{P}}\right)}\vol(P)= \frac{1}{2^n n!},\quad P\in\mathcal{C}^p\left(\R^{2n}\right).\label{eq:proofmain22}\eeq

Moreover, Theorem \ref{Thm:MainProperty1} guarantees the following: $P^*$ is a solution of \eqref{eq:proofmain21} if and only if $P^*$ is a solution of \eqref{eq:proofmain22}.
\epf

\section{Proof of Theorem \ref{Thm:main3}}\label{Sec:main3}

In order to prove Theorem \ref{Thm:main3}, we need the following propositon:

\bprop\label{Prop:characteristicnotranslate}
Let $C\subset\R^{2n}$ be a convex body and $x$ a closed characteristic on $\partial C$. Then, $x$ cannot be translated into $\mathring{C}$.
\eprop

\bpf
Let us assume that $x$ can be translated into $\mathring{C}$. Let $\widetilde{T}>0$ be the period of $x$. Because of
\beqq \dot{x}(t)\in J\partial H_C(x(t))\quad \text{a.e. on }[0,\widetilde{T}],\eeqq
there is a vector-valued function $n_C$ on $\partial C$ such that
\beqq \dot{x}(t)=Jn_C(x(t))\quad \text{a.e. on }[0,\widetilde{T}]\eeqq
with
\beqq n_C(x(t))\in \partial H_C(x(t))\eeqq
for all $t\in[0,\widetilde{T}]$ for which $\dot{x}(t)$ exists and
\beqq n_C(x(t))=0\eeqq
for all $t\in[0,\widetilde{T}]$ for which $\dot{x}(t)$ does not exist.

Then, the convex cone $U$ spanned by
\beqq n_C(x(t))\in N_C(x(t)),\;t\in [0,\widetilde{T}],\eeqq
has the property
\beqq \forall u\in U\setminus\{0\}:\;-u\notin U,\eeqq
since, otherwise, one could find points on $x$ and $C$-supporting hyperplanes through these points with the property that the intersection of the $C$-containing half-spaces bounded by these hyperplanes is nearly bounded (what would imply that $x$ cannot be translated into $\mathring{C}$). By the convexity of $U$, this implies that
\beqq \int_0^{\widetilde{T}} n_C(x(t))\, \mathrm{d}t \neq 0,\eeqq
and therefore
\beqq \int_0^{\widetilde{T}} Jn_C(x(t))\, \mathrm{d}t \neq 0.\eeqq
Since $x$ is a closed characteristic on $\partial C$, $x$ fulfills $x(0)=x(\widetilde{T})$. This implies
\beqq 0=x(\widetilde{T})-x(0)=\int_0^{\widetilde{T}} \dot{x}(t)\, \mathrm{d}t = \int_0^{\widetilde{T}} Jn_C(x(t))\, \mathrm{d}t \neq 0,\eeqq
a contradiction. Therefore, it follows that $x$ cannot be translated into $\mathring{C}$. 
\epf

We now consider the operator norm of the complex structure/symplectic matrix $J$. It is given by:
\beqq ||J||_{C^\circ\rightarrow C}=\sup_{||v||_{C^\circ}\leq 1}||Jv||_{C}=\sup_{\mu_{C^\circ}(v)\leq 1} \mu_C(Jv).\eeqq
We derive the following lemma:

\blem\label{Lem:Karasev}
Let $C\subset\R^{2n}$ be a convex body and $x$ a closed characteristic on $\partial C$ which has period $\widetilde{T}>0$. Then, we have
\beqq \mu_C(\dot{x}(t)) \leq ||J||_{C^\circ \rightarrow C}\quad \text{a.e. on }[0,\widetilde{T}].\eeqq
\elem

\bpf
Since $x$ is a closed characteristic on $\partial C$, we have
\beqq \dot{x}(t)\in J\partial H_C(x(t))\quad \text{a.e. on }[0,\widetilde{T}].\eeqq
This implies
\beqq H_{C^\circ}(-J\dot{x}(t))\in H_{C^\circ}(\partial H_C(x(t)))\quad \text{a.e. on }[0,\widetilde{T}].\eeqq
Using Proposition \ref{Prop:ConsEuler}, we conclude
\beqq H_{C^\circ}(-J\dot{x}(t))=H_C(x(t))=\frac{1}{2}\quad \text{a.e. on }[0,\widetilde{T}],\eeqq
i.e.,
\beqq \mu_{C^\circ}(-J\dot{x}(t))=1\quad \text{a.e. on }[0,\widetilde{T}].\eeqq
Therefore, for
\beqq v(t):=-J\dot{x}(t) \quad \text{a.e. on }[0,\widetilde{T}],\eeqq
we have
\beqq \mu_{C^\circ}(v(t))=1 \; \text{ and } \; Jv(t)=\dot{x}(t) \quad \text{a.e. on }[0,\widetilde{T}]\eeqq
and consequently
\beqq \mu_C(\dot{x}(t)) \leq \sup_{\mu_{C^\circ}(v)\leq 1}\mu_C(Jv) = ||J||_{C^\circ\rightarrow C}\quad \text{a.e. on }[0,\widetilde{T}].\eeqq
\epf

\bpf[Proof of Theorem \ref{Thm:main3}]
By Theorem \ref{Thm:relationcapacity}, we have
\beq c_{EHZ}(C\times C^\circ)=\min_{q\in F^{cp}(C)}\ell_{C^\circ}(q).\label{eq:main31}\eeq
Let $x$ be an action-minimizing closed characteristic on $\partial C$, i.e., $x$ fulfills
\beqq \dot{x}\in J\partial H_C(x)\quad \text{a.e.}\eeqq
and minimizes the action with
\beq \A(x)=-\frac{1}{2}\int_0^{\widetilde{T}} \langle J\dot{x}(t),x(t)\rangle\, \mathrm{d}t=\int_0^{\widetilde{T}} H_C(x)\, \mathrm{d}t=\frac{\widetilde{T}}{2},\label{eq:main33}\eeq
where we used Euler's identity (cf.\;Proposition \ref{Prop:Euleridentity}) to derive
\beqq \langle y,x(t)\rangle = H_C(x(t))\quad \forall y \in \partial H_C(x(t)).\eeqq
Then, since $x$ is in $\partial C$ and, by Proposition \ref{Prop:characteristicnotranslate}, cannot be translated into $\mathring{C}$, \eqref{eq:main31} together with
\beq \min_{q\in F^{cp}(C)}\ell_{C^\circ}(q)=\min_{q\in F^{cc}(C)}\ell_{C^\circ}(q)\label{eq:main34}\eeq
(cf.\;Proposition \ref{Prop:selfevident2}) implies that
\beqq c_{EHZ}(C\times C^\circ)\leq \ell_{C^\circ}(x) = \int_0^{\widetilde{T}} \mu_C(\dot{x}(t))\, \mathrm{d}t.\eeqq
Using Lemma \ref{Lem:Karasev} and \eqref{eq:main33}, we conclude
\begin{align*}
c_{EHZ}(C\times C^\circ)\leq \int_0^{\widetilde{T}} \mu_C(\dot{x}(t))\, \mathrm{d}t &\leq \int_0^{\widetilde{T}} ||J||_{C^\circ \rightarrow C}\, \mathrm{d}t\\
& = \widetilde{T} ||J||_{C^\circ \rightarrow C}\\
& = 2 \A(x)||J||_{C^\circ \rightarrow C}\\
& = 2c_{EHZ}(C)||J||_{C^\circ \rightarrow C}.
\end{align*}
This implies
\beqq \widetilde{R}_C=\frac{c_{EHZ}(C)}{c_{EHZ}(C\times C^\circ)} \geq \frac{1}{2||J||_{C^\circ \rightarrow C}}.\eeqq
Therefore, Viterbo's conjecture for convex bodies in $\R^{2n}$ is equivalent to
\beqq \vol(C)\geq \frac{c_{EHZ}(C)^n}{n!}=\frac{\widetilde{R}_C^n c_{EHZ}(C\times C^\circ)^n}{n!},\quad C\in\mathcal{C}\left(\R^{2n}\right).\eeqq
By referring to Proposition \ref{Prop:invViterbo}, we can assume
\beqq \vol(C)=1\eeqq
without loss of generality and get
\beqq c_{EHZ}(C\times C^\circ)\leq \frac{\sqrt[n]{n!}}{\widetilde{R}_C},\quad C\in\mathcal{C}\left(\R^{2n}\right),\eeqq
which, by Theorem \ref{Thm:relationcapacity}, is equivalent to
\beqq \max_{\vol(C)=1}\; \min_{q\in F^{cp}(C)}\ell_{C^\circ}(q) \leq \frac{\sqrt[n]{n!}}{\widetilde{R}_C},\quad C\in\mathcal{C}\left(\R^{2n}\right).\eeqq
By Theorem \ref{Thm:MainProperty1}, this is equivalent to
\beqq \min_{C\in A\left(C^\circ,\frac{\sqrt[n]{n!}}{\widetilde{R}_C}\right)}\vol(T)\geq 1,\quad C\in\mathcal{C}\left(\R^{2n}\right),\eeqq
and, after applying Proposition \ref{Prop:Ahomogenity2}, to
\beqq \min_{C\in A\left(C^\circ,\frac{1}{\widetilde{R}_C}\right)}\vol(T)\geq \frac{1}{n!},\quad C\in\mathcal{C}\left(\R^{2n}\right).\eeqq

By similar reasoning, Theorem \ref{Thm:MainProperty1} also guarantees the equivalence of
\beq \max_{\vol(C)=1}\; \min_{q\in F^{cp}(C)}\ell_{C^\circ}(q) = \frac{\sqrt[n]{n!}}{\widetilde{R}_C},\quad C\in\mathcal{C}\left(\R^{2n}\right),\label{eq:main35}\eeq
and
\beq \min_{C\in A\left(C^\circ,\frac{1}{\widetilde{R}_C}\right)}\vol(T) = \frac{1}{n!},\quad C\in\mathcal{C}\left(\R^{2n}\right).\label{eq:main36}\eeq

Moreover, Theorem \ref{Thm:MainProperty1} guarantees the following: $C^*$ is a solution of \eqref{eq:main35} if and only if $C^*$ is a solution of \eqref{eq:main36}.
\epf

\section{Justification of Conjectures \ref{Conj:Wetzelsproblem} and \ref{Conj:WetzelsproblemBilliard}}\label{Sec:Wetzelproblem}

We transfer Viterbo's conjecture onto Wetzel's problem. For that, we define 
\beqq y:=\min_{K\in A\left(B_1^2,1\right)} \vol(K)\eeqq
and let $K^*\subset\R^2$ be an arbitrarily chosen convex body of volume $y$. Then, applying Theorems \ref{Thm:relationcapacity} and \ref{Thm:MainProperty1}, we have
\begin{align*}
\frac{c_{EHZ}\left(K^*\times B_1^2\right)^2}{2} &\leq \max_{\vol(K)=y}\frac{c_{EHZ}\left(K\times B_1^2\right)^2}{2}\\
& = \max_{\vol(K)=y}\;\min_{q\in F^{cp}(K)}\frac{\ell_{B_1^2}(q)^2}{2}\\
& =\frac{1}{2}.
\end{align*}
Further, we have
\beqq \vol\left(K^* \times B_1^2\right)=\pi y.\eeqq
The truth of Viterbo's conjecture requires
\beqq \vol\left(K^* \times B_1^2\right)\geq \frac{c_{EHZ}\left(K^* \times B_1^2\right)^2}{2},\eeqq
i.e., $\pi y \geq \frac{1}{2}$, which means
\beqq y\geq \frac{1}{2\pi}\approx 0.15915.\eeqq
Theorem \ref{Thm:MainProperty1} also guarantees the sharpness of this estimate.

Together with Theorem \ref{Thm:relationcapacity}, this justifies the formulation of Conjectures \ref{Conj:Wetzelsproblem} and \ref{Conj:WetzelsproblemBilliard}.

\section{Proof of Theorem \ref{Thm:analogueFinchWetzel}}\label{Sec:escape}

In order to prove Theorem \ref{Thm:analogueFinchWetzel}, we start with the two following obvious observations:

\bprop\label{Prop:selfevident}
Let $K\subset\R^n$ be a convex body. Then we have
\beqq \{\text{closed Minkowski escape paths for $K$}\} = F^{cc}(K).\eeqq
\eprop

\bpf
The statement follows directly by recalling that a closed Minkowski escape path is a closed curve whose all translates intersect $\partial K$ and therefore, equivalently, cannot be translated into $\mathring{K}$.
\epf

\bprop\label{Prop:selfevident2}
Let $K,T\subset\R^n$ be convex bodies. Then we have
\beqq \min_{q\in F^{cc}(K)}\ell_T(q) = \min_{q\in F^{cp}(K)}\ell_T(q).\eeqq 
\eprop

\bpf
Since
\beqq F^{cp}(K) \subset F^{cc}(K),\eeqq
it suffices to find for every closed curve $q\in F^{cc}(K)$ a closed polygonal curve $\widetilde{q}\in F^{cp}(K)$ with
\beq \ell_T(\widetilde{q}) \leq \ell_T(q).\label{eq:selfevident2}\eeq

If $q$ cannot be translated into $\mathring{K}$, then by the remark beyond \cite[Lemma 2.1]{KruppRudolf2020}, there are $n+1$ points on $q$ that cannot be translated into $\mathring{K}$. By connecting these points, we obtain a closed polygonal curve in $F^{cp}(K)$ which we call $\widetilde{q}$. By the subadditivity of the Minkowski functional, it follows \eqref{eq:selfevident2}.
\epf

Based on these propositions, we can prove the analogue to Theorem \ref{Thm:FinchWetzel}:

\bpf[Proof of Theorem \ref{Thm:analogueFinchWetzel}]
We first use Proposition \ref{Prop:selfevident} in order to reduce the statement of Theorem \ref{Thm:analogueFinchWetzel} to: An/The $\ell_T$-minimizing closed curve in $F^{cc}(K)$ has $\ell_T$-length $\alpha^*$ if and only if $\alpha^*$ is the largest $\alpha$ for which
\beq K\in A(T,\alpha).\label{eq:analogue0}\eeq

First, let us asssume that $\alpha^*$ is the $\ell_T$-length of an/the $\ell_T$-minimizing closed curve in $F^{cc}(K)$. Then, from Proposition \ref{Prop:selfevident2}, we know that there is a closed polygonal curve
\beqq q^*\in F^{cp}(K) \; \text{ with } \; \ell_T(q^*)=\alpha^*,\eeqq
i.e., $q^*$ is a minimizer of
\beqq \min_{q\in F^{cp}(K)} \ell_T(q).\eeqq
Then it follows from Proposition \ref{Prop:p^*minF} that
\beqq K\in A(T,\ell_T(q^*))=A(T,\alpha^*).\eeqq
Let $\alpha > \alpha^*$. If
\beq K\in A(T,\alpha),\label{eq:analogue1}\eeq
then
\beqq L_T(\alpha)\subseteq C(K),\eeqq
i.e., every closed curve of $\ell_T$-length $\alpha$ can be covered by a translate of $K$. This implies that every closed curve of $\ell_T$-length $\lambda \alpha $, $\lambda <1$, can be covered by a translate of $\mathring{K}$. From this we conclude
\beqq q^* \notin F^{cp}(K).\eeqq
Therefore, there is no $\alpha > \alpha^*$ for which \eqref{eq:analogue1} is fulfilled, i.e., $\alpha^*$ is the largest $\alpha$ for which \eqref{eq:analogue0} holds.

Conversely, if $\alpha^*$ is the largest $\alpha$ for which \eqref{eq:analogue0} holds. Then, there is a closed curve $q^*$ with
\beq q^*\in F^{cc}(K)\cap C(K) \; \text{ and } \; \ell_T(q^*)=\alpha^*.\label{eq:analogue2}\eeq

Otherwise, if not, then one has
\beqq q\in C(K)\setminus F^{cc}(K)\eeqq
for all closed curves $q$ of $\ell_T$-length $\alpha^*$. This implies
\beqq q\in C(\mathring{K})\eeqq
for all closed curves $q$ of $\ell_T$-length $\alpha^*$. But then there is a $\lambda >1$ such that
\beqq \lambda q \in C(\mathring{K})\eeqq
for all closed curves of $\ell_T$-length $\alpha^*$. But this is a contradiction to the fact that $\alpha^*$ is the largest $\alpha$ for which \eqref{eq:analogue0} holds.

Now, if
\beqq \min_{q\in F^{cc}(K)}\ell_T(q)=:\widetilde{\alpha}<\alpha^*\eeqq
and $\widetilde{q}$ is a minimizer of the left side, then it follows
\beqq \widetilde{q}\in C(K)\eeqq
because, due to Proposition \ref{Prop:cases1}, with $\widetilde{\alpha}<\alpha^*$ one has
\beqq K\in A(T,\alpha^*)\subseteq A(T,\widetilde{\alpha}).\eeqq
Then, with Lemma \ref{Lem:intersectionFC}, there is a $\lambda >1$ such that
\beqq \lambda \widetilde{q}\in F^{cc}(K)\setminus C(K)\eeqq
with
\beqq \ell_T(\lambda\widetilde{q})<\alpha^*.\eeqq
But this is a contradiction to the fact that every closed curve of $\ell_T$-length $\leq \alpha^*$ can be covered by a translate of $K$. Therefore, it follows
\beqq \min_{q\in F^{cc}(K)}\ell_T(q) \geq \alpha^*,\eeqq
and together with \eqref{eq:analogue2}, we conclude that
\beqq \min_{q\in F^{cc}(K)}\ell_T(q) = \alpha^*.\eeqq
\epf

The proof of Corollary \ref{Cor:analogue} follows immediately:

\bpf[Proof of Corollary \ref{Cor:analogue}]
The proof follows directly by combining Proposition \ref{Prop:selfevident2}, \cite[Theorem 3.12]{KruppRudolf2022}, and Theorem \ref{Thm:analogueFinchWetzel}.
\epf

\section[Approach for improving the lower bound in Wetzel's problem]{Computational approach for improving the lower bound in Wetzel's problem}\label{Sec:OptimizationProblem}

In this section, we aim to present a computational approach for improving the best lower bound in Wetzel's problem, which, as stated in Theorem \ref{Thm:WetzelBezdekConnelly}, is due to Wetzel himself (cf.\;\cite{Wetzel1973}). But not only that, our approach most likely also allows to find, more generally, lower bounds in Minkowski worm problems. By Theorem \ref{Thm:MainProperty1}, these lower bounds eventually translate into upper bounds for systolic Minkowski billiard inequalities as well as for Viterbo's conjecture for convex Lagrangian products.

The main idea of this approach is inspired by a series of works related to the search for area-minimizing convex hulls of closed curves in the plane which are allowed to be translated and rotated. Since the area-minimizing convex cover for a set of closed curves is, equivalently, the area-minimizing convex hull of these closed curves (note that this observation has already used within the proof of Corollary \ref{Cor:optimizationmain1}), these works treat the question of lower bounds for the following version of Moser's worm problem in which \textit{closed} arcs are considered:\\

\par
\begingroup
\leftskip=1cm
\noindent \textit{Find a/the convex set of least area that contains a congruent copy of each closed arc in the plane of length one.}\\
\par
\endgroup

\noindent In \cite{ChakerianKlamkin1971} (applying results from \cite{FaryRedei1950}), the first lower bound for the area was found considering the convex hull of a circle and a line segment. In \cite{FurediWetzel2011}, this lower bound was improved by first considering a circle and a certain rectangle and later a circle and a curvilinear rectangle. The latest improvements are due to Grechuk and Som-am who in \cite{GrechukSomam2019} considered the convex hull of a circle, an equilateral triangle and a certain rectangle, and in \cite{GrechukSomam2020} the convex hull of a circle, a certain rectangle, and a line segment. However, in order to adapt these approaches to our setting, in the details, we have to make some changes.

But let us first start with some underlying considerations (as in the proof of Corollary \ref{Cor:optimizationmain1}) in the most general case: For arbitrary convex body $T\subset\R^n$, we ask for lower bounds of
\beq \min_{K\in A(T,1)}\vol(K).\label{eq:approach1}\eeq
By referring to the above mentioned main idea, we start by noting that for
\beqq q_1,...,q_k\in L_T(1)\eeqq
we have
\beq \min_{(a_1,...,a_k)\in (\R^n)^k}\;\vol\left(\conv\{q_1+a_1,...,q_k+a_k\}\right) \leq \min_{K\in A(T,1)}\vol(K).\label{eq:approach2}\eeq
This estimate can be further improved by
\beqq \max_{q_1,...,q_k\in L_T(1)}\;\min_{(a_1,...,a_k)\in (\R^n)^k}\;\vol(\conv\{q_1+a_1,...,q_k+a_k\}) \leq \min_{K\in A(T,1)}\vol(K),\eeqq
so that, eventually, we get
\beqq \min_{a_q\in\R^n}\vol\bigg(\conv\bigg\{ \bigcup_{q\in L_T(1)} (q+a_q)\bigg\}\bigg)=\min_{K\in A(T,1)}\vol(K),\eeqq
where the minimum on the left runs for every $q\in L_T(1)$ over all possible translations in $\R^n$.

Let us now exemplary show how \eqref{eq:approach2} can be used to calculate lower bounds of \eqref{eq:approach1} within the setting of Wetzel's problem, i.e., $n=2$ and $T=B_1^2$.

Let $q_1$ be the boundary of $B_{\frac{1}{2\pi}}^2$,
\beqq q_2=q_2(t_1,t_2,\theta)\eeqq
the boundary of an equilateral triangle $T_{t_1,t_2,\frac{1}{3},\theta}$ with mass point $(t_1,t_2)$, side length $\frac{1}{3}$, and angle $\theta$ between one of the sides and the horizontal line, and let
\beqq q_3=q_3\left(r_1,r_2,\widehat{q}\right)\eeqq
be the boundary of a rectangle $R_{r_1,r_2,1,\widehat{q}}$ with middle point $(r_1,r_2)$, perimeter $1$, and quotient of the side lengths $\widehat{q}$. 

\begin{figure}[h!]
\centering
\def\svgwidth{380pt}
\begingroup%
  \makeatletter%
  \providecommand\color[2][]{%
    \errmessage{(Inkscape) Color is used for the text in Inkscape, but the package 'color.sty' is not loaded}%
    \renewcommand\color[2][]{}%
  }%
  \providecommand\transparent[1]{%
    \errmessage{(Inkscape) Transparency is used (non-zero) for the text in Inkscape, but the package 'transparent.sty' is not loaded}%
    \renewcommand\transparent[1]{}%
  }%
  \providecommand\rotatebox[2]{#2}%
  \newcommand*\fsize{\dimexpr\f@size pt\relax}%
  \newcommand*\lineheight[1]{\fontsize{\fsize}{#1\fsize}\selectfont}%
  \ifx\svgwidth\undefined%
    \setlength{\unitlength}{285.00248127bp}%
    \ifx\svgscale\undefined%
      \relax%
    \else%
      \setlength{\unitlength}{\unitlength * \real{\svgscale}}%
    \fi%
  \else%
    \setlength{\unitlength}{\svgwidth}%
  \fi%
  \global\let\svgwidth\undefined%
  \global\let\svgscale\undefined%
  \makeatother%
  \begin{picture}(1,0.75494953)%
    \lineheight{1}%
    \setlength\tabcolsep{0pt}%
    \put(0,0){\includegraphics[width=\unitlength,page=1]{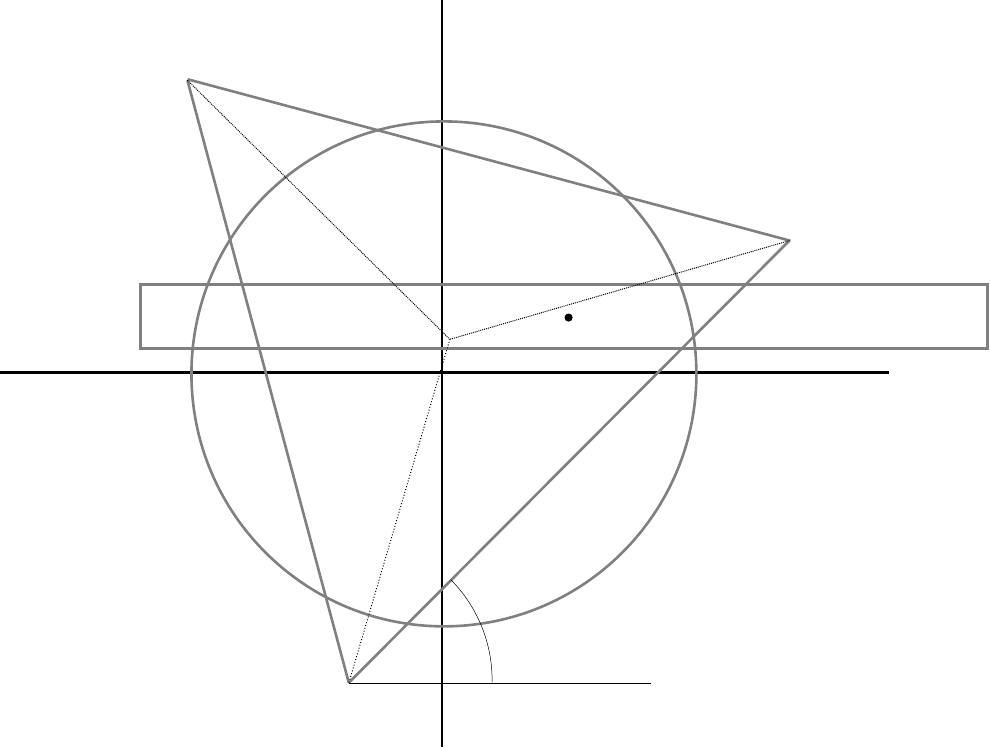}}%
    \put(0.50911268,0.08245986){\color[rgb]{0,0,0}\makebox(0,0)[lt]{\lineheight{1.25}\smash{\begin{tabular}[t]{l}$\theta$\end{tabular}}}}%
    \put(0.58359779,0.42404799){\color[rgb]{0,0,0}\makebox(0,0)[lt]{\lineheight{1.25}\smash{\begin{tabular}[t]{l}$(r_1,r_2)$\end{tabular}}}}%
    \put(0.44775877,0.43989076){\color[rgb]{0,0,0}\makebox(0,0)[lt]{\lineheight{1.25}\smash{\begin{tabular}[t]{l}$(t_1,t_2)$\end{tabular}}}}%
    \put(0,0){\includegraphics[width=\unitlength,page=2]{approach_grey.pdf}}%
    \put(0.8209681,0.43150312){\color[rgb]{0,0,0}\makebox(0,0)[lt]{\lineheight{1.25}\smash{\begin{tabular}[t]{l}$R_{r_1,r_2,1,q}$\end{tabular}}}}%
    \put(0.48638458,0.51890839){\color[rgb]{0,0,0}\makebox(0,0)[lt]{\lineheight{1.25}\smash{\begin{tabular}[t]{l}$T_{t_1,t_2,\frac{1}{3},\theta}$\end{tabular}}}}%
    \put(0.55029375,0.20406154){\color[rgb]{0,0,0}\makebox(0,0)[lt]{\lineheight{1.25}\smash{\begin{tabular}[t]{l}$B_{\frac{1}{2\pi}}^2$\end{tabular}}}}%
    \put(0,0){\includegraphics[width=\unitlength,page=3]{approach_grey.pdf}}%
  \end{picture}%
\endgroup%

\caption{Illustration of the convex hull of $B_{\frac{1}{2\pi}}^2$, $R_{r_1,r_2,1,\widehat{q}}$ and $T_{t_1,t_2,\frac{1}{3},\theta}$.}
\label{img:approach}
\end{figure}

Then, by definition, we have
\beqq q_1,\,q_2(t_1,t_2,\theta),\,q_3\left(r_1,r_2,\widehat{q}\right)\in L_{B_1^2}(1)\eeqq
for all
\beqq t_1,t_2\in\R, \, \theta \in \left[0,\frac{3\pi}{4}\right],\,r_1,r_2\geq 0,\,\widehat{q} >0\eeqq
and \eqref{eq:approach2} (because of $\theta\in\left[0,\frac{3\pi}{4}\right]$ and $\widehat{q}>0$, one has $k=\infty$) becomes
\begin{align*}
\max_{\theta\in[0,\frac{3\pi}{4}],\;\widehat{q}>0}\;\min_{t_1,t_2\in\R, \; r_1,r_2\geq 0}\; \vol &\left(\conv\left\{B_{\frac{1}{2\pi}}^2,T_{t_1,t_2,\frac{1}{3},\theta},R_{r_1,r_2,1,\widehat{q}}\right\}\right)\\
&\leq \min_{K\in A(B_1^2,1)}\vol(K).
\end{align*}
Then, one can define
\beqq f\left(t_1,t_2,r_1,r_2,\theta,\widehat{q}\right):=\vol \left(\conv\left\{B_{\frac{1}{2\pi}}^2,T_{t_1,t_2,\frac{1}{3},\theta},R_{r_1,r_2,1,\widehat{q}}\right\}\right)\eeqq
which is a convex function with respect to the first four coordinates $(t_1,t_2,r_1,r_2)$ (this can be shown similar to in \cite{GrechukSomam2019}) and compute
\beqq \max_{\theta\in[0,\frac{3\pi}{4}],\;\widehat{q}>0}\;\min_{t_1,t_2\in\R, \; r_1,r_2\geq 0}\; f\left(t_1,t_2,r_1,r_2,\theta,\widehat{q}\right).\eeqq

We leave it at that, starting with \eqref{eq:approach2}, gives us the ability to tackle many different Minkowski worm problems--in any dimension, for many different $T$s and by using diverse closed curves
\beqq q_1,...,q_k\in L_T(1).\eeqq

\section*{Acknowledgement}
This research is supported by the SFB/TRR 191 'Symplectic Structures in Geometry, Algebra and Dynamics', funded by the \underline{German Research Foundation}, and was carried out under the supervision of Alberto Abbondandolo (Ruhr-Universit\"at Bochum). The author is thankful to the supervisor's support.


\medskip


\section*{Daniel Rudolf, Ruhr-Universit\"at Bochum, Fakult\"at f\"ur Mathematik, Universit\"atsstrasse 150, D-44801 Bochum, Germany.}
\center{E-mail address: daniel.rudolf@ruhr-uni-bochum.de}

\end{document}